\documentstyle[12pt]{article}
\newtheorem{ex}{Example}[section]
\newtheorem{prob}{Problem}
\newtheorem{teo}{Theorem}
\newtheorem{pr}{Proposition}[section]

\newtheorem{lm}{Lemma}[section]
\newtheorem{rem}{Remark}[section]
\newtheorem{remo}{Remark}

\newcommand{\gr}{{\mbox{gr}\,}}

\newcommand{\Hom}{{\mbox{Hom}\,}}

\begin{document}
\title{Invariants of mixed representations of quivers II: defining relations
and applications}
\author{A.N.Zubkov\\
644099, Omsk-99, Tuhachevskogo embankment 14,\\
Omsk State Pedagogical University ,\\
chair of geometry, e-mail: zubkov@iitam.omsk.net.ru}
\date{}
\maketitle

\section*{Introduction}

In this article we give the complete answer for the problem 1 statetd in the
\cite{zub7}. We recall necessary definitions and notations
(see also \cite{gab, don1, prb1, prb2, zub4}).

A {\it quiver} is a
quadruple $Q=(V,A,h,t)$, where $V$ is a vertex set, $A$ is an arrow set of
$Q$, and the maps $i,t:A\to V$ associate to each arrow $a\in A$ its
origin $i(a)\in V$ and its end $t(a)\in V$. We enumerate elements of the
vertex set as $V=\{1,\ldots ,n\}$.

We consider a collection of vector
spaces $E_1,\ldots ,E_n$ over an algebraically closed field $K$. Let $\dim
E_1=d_1,\ldots ,\dim E_n=d_n$. Denote by ${\bf d}$ the vector $(d_1,\ldots
,d_n)$. This vector is called a {\it dimension vector}.
For two dimension vectors ${\bf d}(1),
{\bf d}(2)$ we write
${\bf d}(1)\geq{\bf d}(2)$ iff $\forall i\in V, d(1)_i\geq d(2)_i$.
Denote by $GL({\bf d})$
the group $GL(E_1)\times \ldots \times GL(E_n)=GL(d_1)\times \ldots \times
GL(d_n)$. The representation space of the quiver $Q$ of dimension ${\bf d}$ is 
$R(Q, {\bf d})=\prod_{a\in A}\Hom_K(E_{i(a)},E_{t(a)})$. The group 
$GL({\bf d})$ acts on $R(Q,{\bf d})$ by the rule:

$$(y_a)_{a\in A}^g=(g_{t(a)}y_ag_{i(a)}^{-1})_{a\in A},
g=(g_1,\ldots , g_n)\in GL({\bf d}),
(y_a)_{a\in A}\in R(Q, {\bf d}).
$$

\noindent For example, if our quiver $Q$ has one vertex and $m$ loops which are
necessarily incident to
this vertex then the ${\bf d}=(d)$-representations space of this quiver is
isomorphic
to the space of $m$ $d\times d$-matrices with respect to the diagonal action
of the group $GL(d)$ by conjugation.

The coordinate ring of the affine variety $R(Q, {\bf d})$
is isomorphic to $K[y_{ij}(a)\mid 1\leq j\leq d_{i(a)}, 1\leq i\leq
d_{t(a)}, a\in A]$. For any $a\in A$ denote by $Y_{{\bf d}}(a)$ the general
matrix $(y_{ij}(a))_{1\leq j\leq d_{i(a)}, 1\leq i\leq d_{t(a)}}$. The
action of $GL({\bf d})$ on $R(Q, {\bf d})$ induces the action on the
coordinate ring by the rule $Y_{{\bf d}}(a)\mapsto g_{t(a)}^{-1}
Y_{{\bf d}}(a)g_{i(a)}, a\in A$. We omit the lower index ${\bf d}$ if it does
not lead to confusion. For example, we write just $Y(a)$ instead of
$Y_{{\bf d}}(a)$.

In \cite{zub7} the concept of a representation of a quiver was generalized as
follows. We partition the vertex set of a given quiver $Q$ into
several disjoint subsets. To be precise, let  $V=
V_{ord}\bigsqcup(\bigsqcup_{q\in\Omega}V_q)$. The vertices from $V_{ord}$
are said to be {\it ordinary}. Moreover, all subsets $V_q$ have cardinality 
two, that is for any $q\in\Omega$ $V_q=\{i_q, j_q\}$.

A dimension vector ${\bf d}$ is said to be {\it compatible} with this
partition of $V$ if for any $q\in\Omega, d_{i_q}=d_{j_q}=d_q$.
The next step is to replace all $E_{j_q}, q\in\Omega$ by their duals. To
indicate that some vertices correspond to the duals of vector spaces we
introduce a new dimension vector
${\bf t} =(t_1,\ldots , t_l)$, where $t_i=d_i$
iff we assign to $i$ the space $E_i$, otherwise $t_i=d_i^*$.
We call ${\bf d}$ the vector underlying ${\bf t}$. These notations
will be used throughout.
By definition, the ${\bf t}$-dimensional representation space of the
quiver $Q$ is equal to the space $R(Q, {\bf t})=\prod_{a\in A}\Hom_K (W_{i(a)}
, W_{t(a)})$, where $W_i=E_i$ iff $t_i=d_i$, otherwise $W_i=E_i^*$.
The space $R(Q, {\bf t})$ is a $G=GL({\bf d})$-module
under the same action:

$$(y_a)_{a\in A}^g=(g_{t(a)}y_ag_{i(a)}^{-1})_{a\in A},
g=(g_1,\ldots , g_l)\in G, (y_a)_{a\in A}\in R(Q, {\bf t}).$$

\noindent If $\Omega=\emptyset$ then ${\bf t}={\bf d}$ and 
$R(Q, {\bf t})=R(Q, {\bf d})$.
Without loss of generality one can identify the coordinate algebras
$K[R(Q, {\bf t})]$ and  $K[R(Q, {\bf d})]$.
Finally, replacing all subfactors $GL(E_{i_q})\times GL(E_{j_q})=GL(d_q)\times
GL(d_q)$ of the group $G=GL({\bf d})$  by
their diagonal subgroups we obtain a new group $H({\bf t})$.
The space $R(Q, {\bf t})$ with respect to the action of  the group
$H({\bf t})$ is called a {\it mixed} representation space of the quiver $Q$
of dimension ${\bf t}$ relative to the partition
$V=V_{ord}\bigsqcup(\bigsqcup_{q\in\Omega}V_q)$.
In \cite{zub7} author formulated the following:

\begin{prob}
What are the generators and the defining relations of the algebra
$J(Q, {\bf t})=K[R(Q, {\bf t})]^{H({\bf t})}$?
\end{prob}

In \cite{zub7} some necessary definitions and 
notations were introduced as well as some auxiliary results
were proved. We remind them since they are
necessary to understand this article.
We start with the notion of the {\it doubled} quiver $Q^{(d)}$. This quiver 
is constructed
with respect to the partition of the vertex set of $Q$ into ordinary 
vertices and
couples $\{i_q, j_q\}, q\in\Omega$. More precisely, 
the vertex set $V^{(d)}$ of $Q^{(d)}$ is equal to $V\bigsqcup V^*_{ord}$, 
where $V^*_{ord}=\{i^*\mid i\in V_{ord}\}$.
Respectively,
the arrow set $A^{(d)}$ of $Q^{(d)}$
is equal to $A\bigsqcup\overline{A}$, where $\overline{A}=\{\bar{a}\mid
a\in A\}$.

Further, if $i(a), t(a)\in V_{ord}$ then $i(\bar{a})=t(a)^*,
t(\bar{a})=i(a)^*$ but if $i(a)$ or $t(a)$ lies in some $V_q, q\in\Omega,$ then

$$i(\bar{a})=\left\{\begin{array}{c}
j_q, t(a)=i_q ,\\
i_q, t(a)=j_q ,
\end{array} \right.
$$

\noindent and symmetrically

$$t(\bar{a})=\left\{\begin{array}{c}
j_q, i(a)=i_q,\\
i_q, i(a)=j_q .
\end{array} \right.
$$

Finally, for any $a\in A^{(d)}$ we suppose $Z(a)=Y(a)$ if $a\in
A$, otherwise $a=\bar{b}, b\in A,$ and $Z(a)=\overline{Y(b)}$, where 
$\overline{Y(b)}$ is the transpose of $Y(b)$.

A product $Z(a_r)\ldots Z(a_1)$ is said to be {\it admissible} if $a_r
\ldots a_1$ is a closed path in $Q^{(d)}$, that is if $t(a_i)=i(a_{i+1}),
i=1, \ldots, m-1,$ and $i(a_1)=t(a_m)$.

The main result of \cite{zub7} is

\begin{teo}
The algebra $J(Q, {\bf t})$ is generated by the elements $\sigma_j(p)$, where
$p$ is an admissible product,  
$\sigma_j$ is $j$-th coefficient of characteristic polynomial, 
$1\leq j\leq\max\limits_{1\leq i\leq n}\{d_i\}$.
\end{teo}

Moreover, it was proved that for any ${\bf d}(1)\geq{\bf d}(2)$ there is
a natural epimorphism $\phi_{{\bf t}(1), {\bf t}(2)} : 
J(Q, {\bf t}(1))\to J(Q, {\bf t}(2))$ which is an
isomorphism on homogeneous components of fixed degree whenever ${\bf t}(1),
{\bf t}(2)$ are sufficiently large. 
For example, if this degree is $r$ then
${\bf t}(1), {\bf t}(2)$ are sufficiently large if all coordinates of their
underlying dimension vectors are not less than $r$ \cite{zub7}.
Since all epimorphisms $\phi_{{\bf t}(1), {\bf t}(2)}$ are
homogeneous we have the countable set of spectrums:

$$
\{ J(Q, {\bf t})(r), \phi_{{\bf t}(1), {\bf t}(2)} \mid {\bf d}
(1)\geq {\bf d}(2)\}, r=0,1,2,\ldots$$

\noindent The inverse limit of $r$-th spectrum denote by $J(Q)(r)$. It is 
clear that
$J(Q)=\oplus_{r\geq 0}J(Q)(r)$ can be endowed with a graded algebra structure 
in obvious way.
The algebra $J(Q)$ is said to be a {\it free} invariant algebra of
mixed representations of the quiver $Q$. One can prove that the algebra $J(Q)$
is a polynomial algebra on infinite many variables (see \cite{zub4}).

By definition
there is a natural epimorphism $J(Q)\to J(Q, {\bf t})$. Denote its kernel
by $T(Q, {\bf t})$. It is clear that this ideal is a $T$-ideal with respect to
some specific set of substitutions on {\it formal} variables corresponding to
arrows of $Q^{(d)}$. More precisely, consider the non-unital {\it path} 
algebra $K<Q^{(d)}>$ of the quiver $Q$. This algebra is generated by the 
elements ${\bf Z}(a), a\in A^{(d)}$. 
The defining relations between these generators are ${\bf Z}(a){\bf Z}(b)=0$ 
iff $i(a)
\neq t(a)$, $a, b\in A^{(d)}$. Let $C(Q)=J(Q)<Q^{(d)}>=J(Q)\otimes 
K<Q^{(d)}>$. 

Notice that $C(Q)$ has an $J(Q)$-algebra involution $\iota$
which is uniquely defined by $\iota : {\bf Z}(a)\mapsto {\bf Z}(\bar{a})$.
It is defined correctly since we suppose that 
$\bar{\bar{a}}=a, a\in A$. It is obvious that  $C(Q)$ is generated by 
the elements
${\bf Y}(a), a\in A,$ as a $J(Q)[\iota ]$-algebra.       

Any non-zero monomial $m={\bf Z}(a_k)\ldots {\bf Z}(a_1)$ corresponds to the 
path
$p=a_k\ldots a_1$ in $Q^{(d)}$. Therefore, one can define its 
origin $i(m)=i(p)=
i(a_1)$ and end $t(m)=t(p)=t(a_k)$. An element $f\in C(Q)$ 
is said to be {\it incident} to
$i\in V^{(d)}$ if all monomials belonging to $f$ are closed pathes in $Q^{(d)}$
starting with $i$. In the same way, we say that  $f\in C(Q)$ is  
{\it passing} from $i$ to $j$, $i, j\in V^{(d)}$, 
if all monomials belonging to $f$ are 
passing from $i$ to $j$. 
The substitution ${\bf Y}(a)\mapsto f_a\in C(Q), a\in A,$ is said to be 
{\it admissible} iff all monomials belonging to $f_a$ has non-zero degree and
pass from $i(a)$ to $t(a)$ \cite{zub1, zub4}. 
  
Specializing all variables
${\bf Y}(a)$ to $Y(a)\in J(Q, {\bf t})$ we see that all admissible 
substitutions induce
endomorphisms of this algebra. Moreover, these endomorphisms
are compatible with the epimorphisms $\phi_{{\bf t}(1), {\bf t}(2)}$. In 
particular, all admissible substitutions induce endomorphisms of $J(Q)$ and
all ideals $T(Q, {\bf t})$ are stable under these endomorphisms. Therefore, it
is possible to regard the ideals $T(Q, {\bf t})$ as $T$-ideals and set the
second part of Problem 1 as follows.

\begin{prob}
What are the generators of $T(Q, {\bf t})$ as a $T$-ideal?
\end{prob}  

In this article we give complete solution of Problem 2. To formulate our 
main result one need more definitions and results from \cite{zub7}.
Decompose the arrow set $A$ into three subsets $A_i, i=1, 2, 3$, where

$$A_1=\{a\in A\mid W_{i(a)}=E_{i(a)}, W_{t(a)}=E_{t(a)}\},
A_2=\{a\in A\mid W_{i(a)}=E_{i(a)}, W_{t(a)}=E_{t(a)}^*\},$$ 
$$A_3=\{a\in A\mid W_{i(a)}=E_{i(a)}^*, W_{t(a)}=E_{t(a)}\}.$$ 

\noindent In other words,
$A_1=\{a\in A\mid i(a), t(a)\in V_{ord}\}, A_2=\{a\in A\mid t(a)=j_q, q\in
\Omega\}, A_3=\{a\in A\mid i(a)=j_q, q\in\Omega\}$. Remark that the case 
$W_{i(a)}=E^*_{i(a)}, W_{t(a)}=E^*_{t(a)}$ can be easily eliminated so we do
not consider it at all 
(see \cite{zub7}).

Fix a {\it multidegree} $\bar{r}=(r_a)_{a\in A}$ and denote $\sum_{a\in A}r_a$
by $r$. Denote by
$J(Q, {\bf t})(\bar{r})$ the homogeneous component of the algebra 
$J(Q, {\bf t})$ of degree $r_a$ in $Y(a)$, $a\in A$.
It was proved in \cite{zub7} that $J(Q, {\bf t})(\bar{r})\neq 0$ iff 
$\sum_{a\in A_2}r_a=\sum_{a\in A_3}r_a=s$. As in \cite{zub7} denote by $t=
r-2s$.  

We extend the set of matrix variables $\{Y(a)\mid a\in A\}$
in the following way. Replace each $Y(a)$ by some new set of 
matrices
having the same size as $Y(a)$. The cardinality of this set is equal to $r_a$.
The same procedure is possible for formal variables ${\bf Y}(a), a\in A$.
Simultaneously, we replace each arrow $a$ by $r_a$ new arrows with the same
origin and end as $a$ and set them in one-to-one correspondence with these new 
matrices. So we get a new quiver $\hat{Q}$.
The vertex set of $\hat{Q}$ coincides with $V$ but the arrow set $\hat{A}$
can be different from $A$.

Set any linear order on $A$. Denote this order by usual
symbol $<$. We enumerate arrows of the quiver
$\hat{Q}$ by the integers $1,\ldots , r$. One can assume that for any 
$a\in A$ the
corresponding set of new arrows is enumerated by the integers from the segment 
$[\dot{a}, a]=[\sum_{b<a}r_b+1,\sum_{b\leq a}r_b]$.
We obtain some {\it arrow specialization} $f: [1, r]=\hat{A}
\to A$ by $f(j)=a$ iff $j\in [\dot{a}, a], a\in A$.
The specialization of matrix variables $Y(j)\mapsto Y(a)$ (or
formal variables ${\bf Y}(j)\mapsto {\bf Y}(a)$) iff
$j\in [\dot{a}, a], a\in A,$ denote by the same symbol $f$.

Without loss of generality it can be assumed that $\forall a\in\hat{A}_1,
b\in\hat{A}_2, c\in\hat{A}_3, a < b < c$. Thus it follows that
$\hat{A}_1=[1,t], \hat{A}_2=[t+1,t+s], \hat{A}_3=[t+s+1, r]$ and $f([1, t])=
A_1$, $f([t+1, s+t])=A_2$, $f([s+t+1, r])=A_3$. It is clear that
$i(j)=i$ or $t(j)=i$ iff $i(f(j))=i$ or $t(f(j))=i$ respectively,
$j\in\hat{A}=[1,\ldots , r], i\in V$.

We set

$$T(i)=
\{j\in\hat{A}\mid t(j)=i\}, I(i)=\{j\in\hat{A}\mid i(j)=i\}, i\in V_{ord},$$

$$T(i_q)=\{j\in\hat{A}\mid t(j)=i_q\}, I(i_q)=\{j\in\hat{A}\mid i(j)=i_q\},$$
$$T(j_q)=\{j\in\hat{A}\mid i(j)=j_q\}, I(j_q)=\{j\in\hat{A}\mid t(j)=j_q\},
q\in\Omega.$$

Denote by $L(Q)$ the subset of the group $S_r$ consisting of all permutations
$\sigma$ which satisfy the conditions:
 
\begin{enumerate}

\item $\forall i\in V_{ord}, \sigma((T(i)\bigcap\hat{A}_1)\bigsqcup (T(i)
\bigcap\hat{A}_3 -s))=(I(i)\bigcap\hat{A}_1)\bigsqcup (I(i)\bigcap
\hat{A}_2 +s)$;

\item $\forall q\in\Omega, \sigma((T(i_q)\bigcap\hat{A}_1)\bigsqcup
(T(i_q)\bigcap\hat{A}_3 -s)\bigsqcup T(j_q))=
(I(i_q)\bigcap\hat{A}_1)\bigsqcup (I(i_q)\bigcap
\hat{A}_2 +s)\bigsqcup I(j_q).$

\end{enumerate}

\noindent These conditions will be called {\it admissibility} conditions. 
As in \cite{zub7} denote the right hand side sets of these
equations by ${\cal I}(i), {\cal I}(q)$. Denote by ${\cal T}(i), 
{\cal T}(q)$ the left hand side sets, that is
the arguments of the permutation $\sigma$. 

Any element $z\in J(\hat{Q})
((1^r))$ is a $K$-linear combination of products of traces 
$tr({\bf Z}(a)\ldots {\bf Z}(b))$, where $a\ldots b$ is a closed path 
in $\hat{Q}$ and
$(1^r)=(\underbrace{1,\ldots , 1}_r)$. As in \cite{zub7} we omit the
functional symbol $tr$ in a record of any $z$ if it does not lead to
confusion. The same convention works in $J(Q)$. For example, if 
$z=(a\ldots b)\ldots (c\ldots d)$ then the specialization
$f$ takes $z$ into $f(z)=(f(a)\ldots f(b))\ldots (f(c)\ldots f(d))$.
An equation $(p)=(q)$ for given closed pathes $p, q$ in $Q^{(d)} \ 
(\hat{Q}^{(d)})$ takes a place iff these pathes
are the same up to a cyclic permutation or involution $\iota$. 
We get an equivalence on the set of closed pathes in $Q^{(d)} 
\ (\hat{Q}^{(d)})$. Any equivalence class is said to be a {\it cycle}.
The given cycle $p$ is called {\it primitive} if it is not a proper power 
\cite{don2}.

There is a $K$-linear isomorphism of vector spaces $tr^* : K[L(Q)]\to 
J(\hat{Q})
((1^r))$ which can be defined as follows. For any $\sigma\in L(Q)$ we have
$tr^*(\sigma)=(a\ldots b)\ldots (c\ldots d)$, where all symbols $a,\ldots b,
\ldots, c,\ldots , d$ lie in the set $[1, r]\bigcup [\bar{1}, \bar{r}]$,
$[\bar{1}, \bar{r}]=\iota([1, r])$, and the set $\{a,\ldots b,
\ldots, c,\ldots , d\}$ has cardinality $r$. Notice that this record is
uniquely defined up to the equivalence mentioned above.
To describe the computation of $tr^*(\sigma)$ one has to define the right 
hand side neighbor of any symbol $j\in [1, r]\bigcup [\bar{1}, \bar{r}]$
in a record of $tr^*(\sigma)$. We list all possibilities for $j$ as follows
(see Proposition 2.5 \cite{zub7}).

Let $j\in [1,r]$ then we have

\begin{enumerate}

\item If $j\in \hat{A}_1$ then $(\ldots jk\ldots)$, where

$$k=\left\{\begin{array}{c}
\sigma^{-1}(j), \ \sigma^{-1}(j)\in \hat{A}_1, \\
\sigma^{-1}(j)+s, \ \sigma^{-1}(j)\in \hat{A}_2, \\
\overline{\sigma^{-1}(j)}, \ \sigma^{-1}(j)\in \hat{A}_3 .
\end{array} \right.
$$

\item If $j\in\hat{A}_2$ then $(\ldots jk\ldots )$, where

$$k=\left\{\begin{array}{c}
\sigma^{-1}(j+s), \ \sigma^{-1}(j+s)\in\hat{A}_1, \\
\sigma^{-1}(j+s)+s, \ \sigma^{-1}(j+s)\in\hat{A}_2, \\
\overline{\sigma^{-1}(j+s)}, \ \sigma^{-1}(j+s)\in\hat{A}_3 .
\end{array} \right.
$$

\item If $j\in\hat{A_3}$ then $(\ldots jk\ldots )$, where

$$k=\left\{\begin{array}{c}
\overline{\sigma(j)}, \ \sigma(j)\in\hat{A}_1, \\
\sigma(j), \ \sigma(j)\in\hat{A}_2, \\
\overline{\sigma(j)-s}, \ \sigma(j)\in\hat{A}_3 .
\end{array} \right.
$$

\end{enumerate}

\noindent If $j=\bar{l}$ then we have the following rules:

\begin{enumerate}

\item If $l\in \hat{A}_1$ then $(\ldots jk\ldots)$, where

$$k=\left\{\begin{array}{c}
\overline{\sigma(l)}, \ \sigma(l)\in \hat{A}_1, \\
\sigma(l), \ \sigma(l)\in \hat{A}_2, \\
\overline{\sigma(l)-s}, \ \sigma(l)\in \hat{A}_3 .
\end{array} \right.
$$

\item If $l\in\hat{A}_2$ then $(\ldots jk\ldots )$, where

$$k=\left\{\begin{array}{c}
\sigma^{-1}(l), \ \sigma^{-1}(l)\in\hat{A}_1, \\
\sigma^{-1}(l)+s, \ \sigma^{-1}(l)\in\hat{A}_2, \\
\overline{\sigma^{-1}(l)}, \ \sigma^{-1}(l)\in\hat{A}_3 .
\end{array} \right.
$$

\item If $l\in\hat{A_3}$ then $(\ldots jk\ldots )$, where

$$k=\left\{\begin{array}{c}
\overline{\sigma(l-s)}, \ \sigma(l-s)\in\hat{A}_1, \\
\sigma(l-s), \ \sigma(l-s)\in\hat{A}_2, \\
\overline{\sigma(l-s)-s}, \ \sigma(l-s)\in\hat{A}_3 .
\end{array} \right.
$$

\end{enumerate}
    
\noindent Following \cite{zub7} we call these rules {\it contracting}.
For other way to compute $tr^*(\sigma)$ see Lemma 3.5 from \cite{zub7}.
Denote $f(tr^*(\sigma))$ by $tr^*(\sigma , f)$. 

Consider the case when $\Omega=\{q\}$ and $V_{ord}=\{i_q\}$.
It is clear that $A_1$ is the set of loops necessary incident to $i_q$.
We suppose additionally that $\mid A_1\mid=\mid A_2\mid=\mid A_3\mid=1$
and ${\bf X}, {\bf Y}, {\bf Z}$ correspond to the single arrows from 
$A_1, A_2$ and $A_3$
respectively. For given $r, s, 2s\leq r,$ we define 
the element 

$$\sigma_{r, s}({\bf X}, {\bf Y}, {\bf Z})=
\frac{1}{t! (s!)^2}\sum_{\sigma\in S_r}
(-1)^{\sigma}tr^*(\sigma, f)$$

\noindent Here $f([1, t])={\bf X}, f([t+1, t+s])={\bf Y}, 
f([t+s+1, r])={\bf Z}$. Notice that in the case
$s=0$ this element coincides with 
$\sigma_r({\bf X})$. Finally, we formulate the main 
result of this article.

\begin{teo}
For any dimension ${\bf t}$ the ideal $T(Q, {\bf t})$ is generated as a
$T$-ideal by the elements $\sigma_r(f), \sigma_{r, s}(f_1, f_2, f_3)$, where
$f, f_1, f_2, f_3\in C(Q)$, $f$ is incident to some $i\in V$ and
$r > d_i$, $f_1$ is incident to $i\in V_{ord}$ or to $i_q, q\in\Omega,$ and
$f_2 (f_3)$ are passing from $i$ or $i_q$ (from $i^*$ or $j_q$) to
$i^*$ or $j_q$ (to $i$ or $i_q$) correspondingly. Moreover, $r > d_i$
and $r > d_q$ respectively.   
\end{teo}

\begin{remo}
If ${\bf t}={\bf d}$ then we get the main result of \cite{zub4}. Notice that
this result was formulated in \cite{zub4} incorrectly. Indeed, it was claimed 
that $f$ in a relation $\sigma_r(f)$ is a monomial incident to $i$. 
It is true if there is a loop $p$ incident to $i$ and then one can replace 
all $\sigma_r(f)$ by $\sigma_r(p)$ up to the obvious substitutions. The same 
remark is for concomitatnts.  
\end{remo}

In the last section we give some applications for orthogonal and symplectic 
invariants of several matrices. More precisely, we describe some approach to 
the problem of computation of defining relations. Notice that in
the characteristic zero case it has been done in \cite{pr}.
I hope to get a complete solution of this problem in the next article.

\section{Preliminaries}

\subsection{Specializations}

For given ${\bf d}(1)\geq {\bf d}(2)$ define
an epimorphism

$$p_{{\bf t}(1), {\bf t}(2)}: K[R(Q, {\bf t}(1))]\to K[R(Q,
{\bf t}(2))]$$

\noindent by the following rule.
Take any arrow $a\in A$. Let $i(a)=i$ and $t(a)=j$. For the sake of simplicity
denote $d_{i}(f)$ and $d_j(f)$ by $m_f$ and $l_f$ respectively, $f=1,2$.
We know that $m_1\geq m_2$ and $l_1\geq l_2$. If $m_f\neq l_f, f=1,2,$ then 
our epimorphism takes $y_{sr}(a)$ to zero iff either $s > l_2$ or $r >
m_2$. On the remaining variables our epimorphism is the identical map.
If $m_f=l_f, f=1,2$ then one can define our epimorphism on the coefficients of
$Y_{{\bf t}(1)}(a)$ in the other way
taking $y_{ss}(a)$ to unit iff $s >m_2=l_2$. We admit both
ways and say that $p_{{\bf t}(1), {\bf t}(2)}$ is {\it standard} if
the second way does not happen at all. Otherwise, $p_{{\bf t}(1), {\bf t}(2)}$
is called {\it non-standard} on arrow $a$ if we define $p_{{\bf t}(1), 
{\bf t}(2)}$ on $Y(a)$ by the second 
way. For example, the epimorphisms
$\phi_{{\bf t}(1), {\bf t}(2)}$ from the introduction are just restrictions of
the standard epimorphisms $p_{{\bf t}(1), {\bf t}(2)}$.   

On the other hand, one can define an isomorphism $i_{{\bf t}(2), {\bf t}(1)}$
of the variety $R(Q, {\bf t}(2))$ onto a closed subvariety of
$R(Q, {\bf t}(1))$ by the dual rule, that is the epimorphism defined above is
the comorphism $i_{{\bf t}(2), {\bf t}(1)}^*$.

By the same way one can define the
isomorphism $j_{{\bf t}(2), {\bf t}(1)}$ of $H({\bf t}(2))$ onto a
closed subgroup 
$H({\bf t}(1))$ just bordering any invertible $d_i(2)\times d_i(2)$
matrix by the $d_i(1)-d_i(2)$ additional rows and columns which are zero
outside of the diagonal tail of length $d_i(1)-d_i(2)$. The entries on this
diagonal tail must all be 1's.

It is not hard to check that
$i_{{\bf t}(2), {\bf t}(1)}(\phi^g)=
i_{{\bf t}(2), {\bf t}(1)}(\phi)^{j_{{\bf t}(2), {\bf t}(1)}(g)}$
for any $g\in H({\bf t}(2))$ and
$\phi\in R(Q, {\bf t}(2))$. The analogous equation is valid for the epimorphism
$p_{{\bf t}(1), {\bf t}(2)}$.

\subsection{Young subgroups}

Decompose the interval $[1, k]=\{1,\ldots , k\}$ into some disjoint
subsets, say
$[1, k]=\bigsqcup_{1\leq j\leq m}T_j$. Define the {\it Young subgroup}
$S_{{\cal T}}=S_{T_1}\times\ldots\times S_{T_m}$ of $S_k$
as the subgroup consisting of all
permutations $\sigma\in S_k$ such that $\sigma(T_j)=T_j, 1\leq j\leq m$.
By definition, $S_T=\{\sigma\in S_k\mid\sigma(T)=T, \forall j\not\in T \
\sigma(j)=j\}$ for arbitrary subset $T$. The subsets $T_1,\ldots , T_m$ are 
said to be the {\it layers} of $S_{{\cal T}}$ \cite{zub1, zub4}.

The group $S_{{\cal T}}$ can be defined in
other way. In fact, let $f$ be a map from $[1,k]$ onto $[1, m]$ defined by
the rule $f(T_j)=j, j=1,\ldots , m$. Then $S_{{\cal T}}=
\{\sigma\in S_k\mid f\circ\sigma=f\}$. Sometimes we will denote the group
$S_{{\cal T}}$ by $S_f$.

For any group $G$ and its subgroup $H$ we denote by $G/H$ some fixed
representative set of the left $H$ cosets if it does not lead to confusion.

For a given superpartition
$\bar{\lambda}=(\lambda_1,
\ldots , \lambda_n)$ which is composed from ordinary (not necessary ordered) 
partitions 
$\lambda_i=(\lambda_{i1},\ldots ,\lambda_{i, s_i}),
i=1,\ldots , n$, denote by
$S_{\bar{\lambda}}$ the Young subgroup of $S_{\mid\bar{\lambda}\mid}$
corresponding to the decomposition of $[1, \mid\bar{\lambda}\mid]$ into
sequential subintervals of lengths $\lambda_{11},\ldots , \lambda_{1, s_1},
\lambda_{21},\ldots ,\lambda_{2, s_2}, \ldots$. 

For example, fix some multidegree $\bar{r}$ and consider the 
superpartition $\Theta=(\lambda_{A_1}, \mu_{A_2}, \\
\gamma_{A_3})$,
where $\lambda_{A_1}=(\lambda_a)_{a\in A_1}, \mu_{A_2}=(\mu_a)_{a\in A_2},
\gamma_{A_3}=(\gamma_a)_{a\in A_3}$, $\lambda_a, \mu_b, 
\gamma_c$ are partitions such that
$\forall a\in A_1, \forall b\in A_2, \forall c\in A_3, \mid\lambda_a \mid=r_a,
\mid\mu_b\mid=r_b, \mid\gamma_c\mid=r_c$. If $f$ is the specialization from
the introduction then $S=S_{\Theta}\leq S_f$ \cite{zub7}. Another important
Young subgroup of $S_r$ is $S_0=S_{\hat{A}_1}\times S_{\hat{A}_2}
\times S_{\hat{A}_3}$. Any $\pi\in S_0$ can be decomposed as $\pi=\pi_1\pi_2
\pi_3$, where $\pi_i\in S_{\hat{A}_i}, i=1,2,3,$. Moreover, any $\pi\in
S_{\hat{A}_i}, i=2,3,$ has a double $\pi^x$, where $x=\prod_{i\in\hat{A}_2}(i
\ i+s)$. As in \cite{zub7} we denote $\pi^x$ by $\pi+s$ if $\pi\in
S_{\hat{A}_2}$, otherwise by $\pi -s$. 
In \cite{zub7} two {\it shift} homomorphisms
$\rho_i : S_0\to S_r, i=1,2,$ were defined by the rule 
$\rho_1(\pi)=\pi_1\pi_2
(\pi_2+s), \rho_2(\pi)=\pi_1(\pi_3 -s)\pi_3, \pi\in S_0$.
The sets ${\cal T}(x) \ ({\cal I}(x)), x\in V_{ord}\bigcup\Omega,$ form a
decomposition of the interval $[1, r]$. Therefore, we have two Young
subgroups $S_{{\cal T}}, S_{{\cal I}}$. It is not hard to prove that
$S_f=\rho_1^{-1}(S_{{\cal I}})\bigcap\rho_2^{-1}(S_{{\cal T}})$.

Divide each
${\cal T}(z)$, $z\in V_{ord}\sqcup\Omega$,  into some sublayers in a 
{\it monotonic} way. In other words, let
${\cal T}(z)=\sqcup _{1\leq j\leq l_z}\beta _{zj}$, where $\max 
\beta_{zj_1}<\min\beta _{zj_2}$ as soon as $j_1<j_2$, and $\max (\min)
\beta_{zj}$ means the maximal (minimal) integer from this sublayer.
Joining over all indices $z$ we obtain a decomposition of the segment
$[1, r]$. Denote by $S_{\bar\beta }$ the Young subgroup corresponding to this 
decomposition. As in \cite{zub1, zub7} this subgroup will be called the 
{\it base} group.
We say that $S_{\bar{\beta}}$ is sufficiently {\it large} for the dimension 
${\bf t}$ if there is some $\beta_{ij}$ or $\beta_{qj}$ such that $\mid
\beta_{ij}\mid > d_i$ or $\mid\beta_{qj}\mid > d_q$ respectively, $i\in 
V_{ord}, q\in\Omega$. 

\subsection{Suitable generators}

For a given base group $S_{\bar{\beta}}$ and elements $\sigma_1\in L(Q), 
\sigma_2\in S_{{\cal T}}$ one can define the so-called {\it suitable} generator

$$z=\frac{1}{\mid S_f\mid}tr^*(\sum_{\tau\in S_{\bar{\beta}}}\sum_{\pi\in S/
(\rho_1^{-1}(S_{\bar{\beta}}^{\sigma_1})\bigcap\rho_2^{-1}
(S_{\bar{\beta}}^{\sigma_2})\bigcap S)}(-1)^{\tau}\rho_1(\pi)
\sigma_1\tau\sigma_2^{-1}\rho_2(\pi)^{-1}, f).$$

It was proved in \cite{zub7} that the ideal
$T(Q, {\bf t})$ is generated as a vector space by those suitable generators 
which belong to
sufficiently large base groups $S_{\bar{\beta}}$ for ${\bf t}$.
Without loss of generality one can suppose that $S=S_f=
\rho_1^{-1}(S_{\bar{\beta}}^{\sigma_1})\bigcap\rho_2^{-1}
(S_{\bar{\beta}}^{\sigma_2})$ and $\sigma_2=1$ \cite{zub7}.
This assumption is correct because $S_f\geq\rho_1^{-1}(S_{\bar{\beta}}
^{\sigma_1})\bigcap\rho_2^{-1}(S_{\bar{\beta}}^{\sigma_2})$ for any base group
$S_{\bar{\beta}}$ and $S_{\bar{\beta}}^{\sigma_2}\leq S_{{\cal T}}$.
If we replace $S_{\bar{\beta}}$ by $S_{\bar{\beta}}^{\sigma_2}$
it means that we make possible any partititons of the sets ${\cal T}(z)$ into
layers of a base group.

Sometimes we will denote an element $\frac{1}{\mid S_f\mid}tr^*(D, f)$, where
$D\in K[L(Q)]$, by $c(D), c(D, f)$ or by $c(D, S_f)$.  

\subsection{Good filtration dimension and a cohomological lemma}

Let $G$ be an algebraic group and $V$ some rational $G$-module. 
The module $V$ is said to have {\it good filtration dimension} 
equal to $r$ if for
any module $U$, admitting a good filtration and for any $k > r$,
$H^k(G, V\otimes U)=0$ but there is at least one $U$ admitting a good 
filtration such that $H^r(G, V\otimes U)\neq 0$ (see \cite{don7, don5, parsh} 
for more definitions). 

\begin{lm}(\cite{don7, parsh}) The good filtration dimension (briefly --
$g.f.d.$) satisfies the following properties:

\begin{enumerate}

\item If $g.f.d.(V)\leq n, g.f.d.(U)\leq m$ then $g.f.d.(V\otimes U)\leq
n+m$.

\item If $0\to V\mapsto X_0\to\ldots\to X_n\to 0$ is an exact
sequence of $G$-modules then $g.f.d.(V)\leq\max\{g.f.d(X_i)+i\}$.

\item If $0\to X_n\to\ldots\to X_0\to V\to 0$ is an exact
sequence of $G$-modules then $g.f.d.(V)\leq\max\{g.f.d(X_i)-i\}$.

\end{enumerate}
  
\end{lm}

Suppose that for another algebraic group $H$ and an $H$-module $W$ we have 
a {\it homomorphism of pairs} $(i, p) : (G, V)\to (H, W)$ \cite{weiss}, 
that is a homomorphism $i : H\to G$ and a linear map $p : V\to W$ 
such that $x p(v)=p(i(x)v), x\in H, v\in V$. We call $p$ the {\it linear part}
of $(i, p)$. 
Using Hochschild complexes \cite{jan} it is easy to see that the map $p$ 
induces a map 
of cohomology groups $H^n(G, V)\to H^n(H, W), n=0,1,2, \ldots$. Moreover, 
we get the following 

\begin{lm} If we have a diagram

$$\begin{array}{ccccccccc}
0 & \to & A & \to & B & \to &
C & \to & 0 \\
& & \downarrow & & \downarrow & & \downarrow & &  \\
0 & \to & A' & \to & B' &
\to & C' & \to & 0
\end{array}
$$

\noindent with exact top and bottom rows 
of $G$ and $H$-modules respectively such that the vertical arrows
are homomorphisms of pairs whose linear parts commute with morhisms in
rows, then we obtain a long commutative diagram

$$\begin{array}{ccccccccccc}
\ldots & \to & H^n(G, A) & \to & H^n(G, B) & \to &
H^n(G, C) & \to & H^{n+1}(G, A) &\to &\ldots\\
& & \downarrow & & \downarrow & & \downarrow & &\downarrow & &  \\
\ldots & \to & H^n(H, A') & \to & H^n(H, B') &
\to & H^n(H, C') & \to & H^{n+1}(H, A') &\to &\ldots
\end{array}
$$

\noindent Here the top and bottom sequences are standard long exact 
sequences of  cohomology groups. 
\end{lm}

Proof. The Lemma  is proved by a routine verification in the 
three-dimensional diagram composed from Hochschild complexes of all modules.
For more details see \cite{weiss}, Theorem 2-1-9.
 
\section{Proof of the main theorem}

We have to show that any suitable generator $z=c(\sum_{\tau\in 
S_{\bar{\beta}}}(-1)^{\tau}\sigma\tau , f)$ is a linear combination of
elements from Theorem 2 possibly multiplied by several $\sigma_j(p), p\in
Q^{(d)},$ with integral coefficients. 

The layers of the group $S_f=\rho_1^{-1}(S_{\bar{\beta}}^{\sigma})
\bigcap\rho_2^{-1}(S_{\bar{\beta}})$ are $A_1$-layers $\sigma(\beta_{uv})
\bigcap\beta_{fg}\bigcap\hat{A}_1$, $A_2$-layers $(\sigma(\beta_{uv})
\bigcap\hat{A}_3-s)\bigcap\sigma(\beta_{fg})\bigcap\hat{A}_2$ or
$A_3$-layers $(\beta_{uv}\bigcap\hat{A}_2+s)\bigcap\beta_{fg}\bigcap
\hat{A}_3$.
  
It is clear that $\sigma$ can be choosen such that its any non-trivial cycle
does not contain two integers from the same  layer of $S_{\bar{\beta}}$. 
Therefore,
the $A_1$-layers can be devided into {\it passive} and
{\it active} layers \cite{zub1, zub4}. More precisely, the passive layers are
$(\beta_{uv}\setminus (\bigsqcup_{fg\neq uv}\sigma(\beta_{fg}))
\bigcap\hat{A}_1=\{i\in\beta_{uv}\bigcap\hat{A}_1 \mid\sigma(i)=i\}$
and all other $A_1$-layers are active.

Let $B$ be a union of all passive layers. By the
definition, $B\subseteq\hat{A}_1$. 
Consider an element $\tau\in
S_r$. The cyclic decomposition of $\tau^{-1}$ has the form 
$(B_1i_1\ldots
B_ki_k)\ldots (B_ti_t\ldots \\
B_si_s)$, where $B_1,\ldots , B_s$ are some
fragments of this decomposition consisting of integers from $B$ and all
$i_1,\ldots , i_s$ are not contained in $B$. Notice that some of
these fragments can be empty or coincided with cycles containing them.

\begin{lm}
The element $tr^*(\tau)$ is uniquely defined by the integers $i_1,\ldots ,
i_s$.
\end{lm}

Proof. We have $tr^*(\tau)=(D_1j_1\ldots D_lj_l)\ldots (D_mj_m
\ldots D_tj_t)$, where any $D_u$ coincides either with some $B_v$ or with its
transposed $\iota(B_v)$. More precisely, we consider the blocks 
$i_kB_1 i_1, \ldots , i_{k-1}B_k i_k,\ldots , i_sB_t i_t ,\ldots , 
i_{s-1}B_s i_s$. By the contracting rules one has to replace each 
$l_v B_v j_v$ by $C_v=l'_v B_v j'_v$
, where

$$l'_v=\left\{\begin{array}{c}
l_v, \mbox{if} \ l_v\in\hat{A}_1, \\
\bar{l}_v, \mbox{if} \ l_v\in\hat{A}_2, \\
l_v-s, \mbox{if} \ l_v\in\hat{A}_3 ,
\end{array} \right.
$$

$$j'_v=\left\{\begin{array}{c}
j_v,\mbox{if} \ j_v\in\hat{A}_1, \\
j_v+s, \mbox{if} \ j_v\in\hat{A}_2, \\
\bar{j}_v, \mbox{if} \ j_v\in\hat{A}_3 .
\end{array} \right.
$$

The next step
is to join all these blocks by same ends. For example, if $C_v=
l'_vB_vj'_v, C_w=l'_wB_wj'_w$ then they can be joined 
if  $j'_w=l'_v$ or $j'_v=l'_w$ as follows:
$l'_wB_wj'_wB_vj'_{v}$ or $l'_vB_vj'_vB_wj'_w$ respectively.
If $\bar{j'}_v=j'_w$ or $\bar{l'}_v=l'_w$ then one has to transpose
one of these blocks and join them in the
obvious way. Continuing this process we will get all cycles of $tr^*(\tau)$
step by step like growing crystalles.

\begin{ex}
Let $r=7, s=2, t=3, \tau=(145)(267), B=\{2, 3\}$. Then $\tau^{-1}$ has the 
following blocks: $1(\emptyset)5, 5(\emptyset)4, 4(\emptyset)1 , 
7(\emptyset)6, 6(2)7, 3$. In other words, $B_1=B_2=B_3=B_4=\emptyset,
B_5=2, B_6=3$. We have $C_1=1(\emptyset)7, C_2=\overline{5}(\emptyset)6, C_3=
\overline{4}(\emptyset)1, C_4=5(\emptyset)\overline{6}, C_5=4(2)\overline{7},
C_6=3$.  By Lemma 2.1 we obtain that $tr^*(\tau)=(17\overline{2}\overline{4})
(\overline{5}6)(3)$. It can easily be checked that the contracting rules
give the same result.    
\end{ex}

\begin{rem}
Using the obvious rule $(CiBj)(ij)=(Ci)(Bj)$, where $C,B$ are fragments of 
given cycles,
we see that $tr^*(\sigma\tau)$ can be computed by the following way. 
Deleting all
integers from $B$ in the cyclic decomposition of
$\tau$ we get some permutation $\tilde{\tau}$ which acts on the set
$[1,r]\setminus B$ such that $\tilde{\tau}^{-1}=(i_1\ldots i_k)\ldots (i_t
\ldots i_s)$. Then $tr^*(\sigma\tau)$  is obtained from $tr^*(\sigma
\tilde{\tau})$
by substituting all passive fragments (or their transposed). 
More precisely, let $(\sigma\tilde{\tau})=(j_1\ldots j_m)\ldots 
(j_n\ldots j_s)$. For given fragment, 
say $B_v$, one can find an integer $x$ such that $j_v=j_x$ and define 
$C_v=j'_{x-1}B_vj'_x$. It remains to join all $C_v$ in the way described above.
\end{rem}

From now on we fix some $\beta_{ij}$ or $\beta_{qj}$ such that $\mid
\beta_{ij}\mid > d_i$ or $\mid\beta_{qj}\mid > d_q$ respectively, $i\in 
V_{ord}, q\in\Omega$. We call it a
{\it selected}  layer but all other layers are called {\it ordinary}. 
Slightly abusing our notations denote the selected layer by $\beta_0$. 
We prove our theorem by induction on two parameters $(r, r-t)$, where $t$ is
the number of ordinary layers.
The first step is to eliminate all passive
layers of $S_f$ except those which are contained in the selected layer.

\begin{pr}
Without loss of generality one can assume that the group $S_f$ has not
any passive layers except those are contained in the selected layer.
\end{pr}

We outline the proof and refer for more details to \cite{zub1, zub4}.
Remark 2.1 is used in the following computations without 
additional references.

Let $\beta_{xj}$ be an ordinary layer and $\beta_{xj}\bigcap\hat{A}_1=
\alpha_1\bigcup\ldots\alpha_l\bigcup\alpha_{l+1}$,
where $\alpha_1,\ldots ,
\alpha_l$ are active layers of $S_f$ but $\alpha_{l+1}\neq
\emptyset$ is a passive one.
We say that $\tau\in S_{{\cal T}}$ has the type $\bar{m}=(m_1,\ldots , m_d,
m_d
,m_{d+1})$, where $d=\mid\beta_{xj}\setminus\alpha_{l+1}\mid$ and
$m_{d+1}=\mid\alpha_{l+1}\mid-\sum_{1\leq w\leq d}m_w$, if all passive
fragments belonging to $\alpha_{l+1}$
in the record of $tr^*(\tau)$ have lengths $m_1,\ldots , m_d$ up to
order. Both $\tau$ and $\sigma\tau$ have the same
type. Denote by $I_{\bar{m}}$ the set consisting of all $\tau\in
S_{\bar{\beta}}$ of type $\bar{m}$. The element

$$z=c(\sum_{\tau\in S_{\bar{
\beta}}}(-1)^{\tau}\sigma\tau, f)=\frac{1}{\mid S_f\mid}tr^*(\sum_{\tau\in 
S_{\bar{\beta}}}(-1)^{\tau}\sigma\tau, f)$$

\noindent can be represented as
$\sum_{\bar{m}}c(g_{\bar{m}}, f)$, where $g_{\bar{m}}=\sum_{\tau\in
I_{\bar{m}}}(-1)^{\tau}\sigma\tau$.

Notice that any summand $g_{\bar{m}}$ is $S_f$-invariant with respect to the 
action $\mu^{\pi}=\rho_1(\pi)\mu\rho_2(\pi)^{-1}, \pi\in S_f, \mu\in S_r$, and
if $\pi\in S_{[1,t]}$ then $\mu^{\pi}=\pi\mu\pi^{-1}$.

Let $\beta_{xj}\setminus\alpha_{l+1}= \{v_1,...,v_d\}$.
One can devide each layer
$\alpha_j$ into some sublayers in the following way.
Integers $v_a$ and $v_b$  from $\alpha_t$
belong to the same new sublayer iff
$m_a = m_b$. We obtain a new group $S_{f'}$, where $f'$ corresponds to
this new partition.

Denote by $M(\bar{m})$ the subset of $I_{\bar{m}}$ consisting of  
all elements such that each fragment of length $m_t$ is the left hand side
neighbor of $v_t$ , $1\leq t\leq d$.
The element $g_{\bar{m}}' = \sum_{\tau\in M(\bar{m})}(-1)^{\tau}
\sigma\tau$ is also
$S_{f'}$-invariant and

$$g_{\bar{m}} = \sum_{x\in S_{f} / S_{f'}}
xg_{\bar{m}}'x^{-1}=
\sum_{x\in S_{f} / S_{f'}}\rho_1(x)g_{\bar{m}}'\rho_2(x)^{-1} .$$

\noindent Using Lemma 3.8 from \cite{zub7} we see that $c(g_{\bar{m}}, f)$ is
obtained from $c(g_{\bar{m}}', f')$ with the help of some glueing of 
variables.
Choose in $\beta_{xj}$ some sublayer $\pi$
such that $\pi\subseteq\alpha_{l+1}$ and $\mid\pi\mid=
m_{d+1}$.
Define a new base subgroup $S_{\bar{\beta}'}$, where
$\bar{\beta}'$ coincides with $\bar{\beta}$ outside of
$\beta_{xj}$  but $\beta_{xj}$ is devided
into two sublayers $\pi$ and $\beta_{xj}\setminus\pi$.
It is clear that the new group $S_{f''}=
\rho_1^{-1}(S_{\bar{\beta}'}^{\sigma})\bigcap\rho_2^{-1}(
S_{\bar{\beta}'})\bigcap S_{f'}$ equals $S_{\bar{\beta}'}\bigcap S_{f'}$.
In other words, $S_{f''}$ conicides with $S_{f'}$ outside of
$\alpha_{l+1}$ and $S_{f''}\bigcap
S_{\alpha_{l+1}}=S_{\pi}\times S_{\alpha_{l+1}\setminus
\pi}$. For the sake of convenience we represent the group $S_{f'}$ as
$S_{g}\times S_{\alpha_{l+1}}$.
We have

$$g_{\bar{m}}' =\sum_{x\in S_{f'}/S_{\bar{f''}}}
xg_{\bar{m}}''x^{-1}=
\sum_{x\in S_{\alpha_{l+1}}/(
S_{\pi}\times S_{\alpha_{l+1}\setminus
\pi})}xg_{\bar{m}}''x^{-1}.
$$

\noindent Here $g_{\bar{m}}'' = 
\sum_{\tau\in S_{\bar{\beta}'}\bigcap M(\bar{m})}
(-1)^{\tau}\sigma\tau$. It remains to consider the element
$c(g_{\bar{m}}'', S_{f''})$. We have

$$
c(g_{\bar{m}}'', S_{f''}) =
c(\sum_{\tau\in S_{\bar{\beta}'\setminus\pi}\cap M(\bar{m})}
(-1)^{\tau}\sigma\tau, S_g\times S_{\alpha_{l+1}\setminus\pi})
\times
c(\sum_{\tau\in S_{\pi}}(-1)^{\tau}\tau, S_{\pi}).
$$

Because of $\alpha_{l+1}\subseteq {\cal T}(x)\bigcap\hat{A}_1={\cal T}(x)
\bigcap\sigma({\cal T}(x))\bigcap\hat{A}_1={\cal T}(x)\bigcap {\cal I}(x)
\bigcap
\hat{A}_1$ we see that for all $j\in\alpha_{l+1}, i(j)=t(j)=x \ (i_x)$,
that is the general matrice $Y$ corresponding to this layer is a square 
matrice. In particular,
$c(\sum_{\tau\in S_{\pi}}(-1)^{\tau}\tau,
S_{\pi})=\sigma_{m_{d+1}}({\bf Y})$. 
If $m_{d+1}>0$ then induction on $r$ completes our eliminating process.

Let $m_{d+1}=0$ and $\bar{\chi}=\bar{\beta}\setminus\alpha_{l+1}=
(\ldots , \beta_{xj}\setminus\alpha_{l+1}, \ldots)$.
Fix a collection of fragments $B_1,..., B_d$ which are contained in
$\alpha_{l+1}, \mid B_i\mid = m_i, 1\leq i\leq d$.
Let $S(\bar{m})$ be a set consisting of all $\tau\in M(\bar{m})$ such that
each $B_t$ is the left hand side neighbor of $v_t$, $1\leq t\leq d$.
Then

$$
c(g_{\bar{m}}',S_{f'}) =
\frac{1}{\mid S_{g}\mid}
\frac{1}{\mid S_{\alpha_{l+1}}\mid}
tr^*(\sum_{x\in S_{\alpha_{l+1}}}x(\sum_{\tau\in S(\bar{m})}
(-1)^{\tau}\sigma\tau)x^{-1}, S_{f'}).
$$

It can easily be checked that
$c(g_{\bar{m}}',S_{f'})$ is obtained from

$$
c(\sum_{\tau\in S_{\bar{\chi}}}(-1)^{\tau}\sigma\tau,
S_g)
$$

\noindent by the substitution ${\bf X}_{g(v_t)}\longrightarrow {\bf Y}^{m_t}
{\bf X}_{g(v_t)},
1\leq t\leq d$. Arguing as above, we see that all
products ${\bf Y}^{m_t}{\bf X}_{g(v_t)}$ are defined correctly and
$S_g = \rho_2^{-1}(S_{\bar{\chi}})
\bigcap\rho_1^{-1}( S_{\bar{\chi}}^{\sigma})$.
Since $\mid\alpha_{l+1}\mid>0$ induction on $r$ completes the proof.

We call any preimage $\sigma^{-1}(\sigma(\beta_{uv})\bigcap\beta_{fg}
\bigcap\hat{A}_1)\subseteq\beta_{uv}$ by an $A_1$-{\it prelayer}. An 
$A_2$-{\it prelayer} $\gamma$ is uniquely
defined by  
$\sigma(\gamma)=(\sigma(\beta_{uv})\bigcap\hat{A}_3-s)\bigcap
\sigma(\beta_{fg})\bigcap\hat{A}_2$
or by
$\sigma(\gamma)=\sigma(\beta_{uv})\bigcap(\sigma(\beta_{fg})
\bigcap\hat{A}_2+s)\bigcap\hat{A}_3$. 
Similarly, a $A_3$-{\it prelayer} 
$\gamma$ is uniquely
defined by $\gamma=(\beta_{uv}\bigcap\hat{A}_2+s)\bigcap\beta_{fg}\bigcap
\hat{A}_3$
or by
$\gamma=\beta_{uv}\bigcap(\beta_{fg}\bigcap\hat{A}_3-s)\bigcap\hat{A}_2$. 
A layer containing
given prelayer is called {\it overlayer}. 
We call a prelayer {\it ordinary} if its overlayer is ordinary. 
The following proposition plays crucial role in the simplification of suitable
generators (see \cite{zub1,zub4}). 

\begin{pr} Any ordinary $A_i$-prelayer coincides with its overlayer up to 
induction on the second parameter, $i=1,2,3$. 
In particular, if $\beta_{uv}$ is an
ordinary layer then $\sigma(\beta_{uv})$ belongs to only one set $\hat{A}_1,
\hat{A}_2$ or $\hat{A}_3$. Analogously, $\beta_{uv}$ belongs to only one of
$\hat{A}_1,\hat{A}_2$ or $\hat{A}_3$. Furthemore, if additionally $\beta_{uv}
\subseteq\hat{A}_1$ then it is covered by only one $A_1$-layer of the 
group $S_f$. 
\end{pr}
  
Proof. The proof is similar to \cite{zub1, zub4} up to some specific details 
which we give below. The proof is step by step splitting of all
ordinary layers which do not satisfy at least one condition mentioned in the
proposition. To be precise, for any such layer, say $\beta_{uv}$,
we extract some proper subset $\gamma\subseteq\beta_{uv}$. Next,  
denote by $H$ the group
$S_{\bar{\beta}\setminus\beta{uv}}\times S_{\gamma}\times S_{\beta_{uv}
\setminus\gamma}$.
The induction step consists of proving that $z$ is a sum of elements 
$c(S(x), f)$ or $c(S'(x), f)$, 
where
$$
S(x)=(-1)^x\sum_{\tau\in H}\sum_{\pi\in S_f/(\rho_1^{-1}(H^{\sigma})\bigcap
\rho_2^{-1}(H^x))}(-1)^{\tau}\rho_1(\pi)\sigma\tau x^{-1}\rho_2(\pi)^{-1}
$$
and
$$
S'(x)=(-1)^x\sum_{\tau\in H}\sum_{\pi\in S_f/(\rho_1^{-1}(H^{\sigma x})
\bigcap\rho_2^{-1}(H))}(-1)^{\tau}\rho_1(\pi)\sigma x\tau\rho_2(\pi)^{-1}.
$$
\noindent Here $x$ runs over some subset of  
$S_{\bar{\beta}}/H=S_{\beta_{uv}}/  
(S_{\gamma}\times S_{\beta_{uv}\setminus\gamma})$. For a given $x$ denote
the groups $\rho_1^{-1}(H^{\sigma})\bigcap\rho_2^{-1}(H^x)$ and 
$\rho_1^{-1}(H^{\sigma x})\bigcap\rho_2^{-1}(H)$ by $S_{f_x}$ and $S'_{f_x}$
correspondingly.
It is clear that both $S_{f_x}, S'_{f_x}$
are contained in $S_f$.

The set $\sigma S_{\bar{\beta}}$ is invariant under substitutions
$a\mapsto \rho_1(\pi)a \rho_2(\pi)^{-1}$ for all $\pi\in S_f$ and parities of
both $a$ and $\rho_1(\pi)a \rho_2(\pi)^{-1}$ are the same.
Therefore, all we need is to prove
that in any $S(x) (S'(x))$ there are no repeated summands
and for any two $S(x), S(y) (S'(x), S'(y))$ either their summands are 
same or they have no summands in common \cite{zub1, zub4}. 

If $z$ can be represented as a sum of the 
elements mentioned above, we say that $z$ admits a {\it disjoin} reduction.
Since the number of ordinary layers of the new base group $H$ is increased,
one can apply induction on the second parameter whenever $z$ admits a disjoin
reduction. Notice that if some ordinary 
overlayer coincides with its $A_i$-prelayer, where $i=1,2,3$, then this 
statement remains true even if we split this overlayer into some sublayers.  

$({\bf i})$ Let $\gamma=\sigma^{-1}(\sigma(\beta_{uv})\bigcap\beta_{fg}
\bigcap\hat{A}_1)$
be an ordinary $A_1$-prelayer. We work with sums $S(x)$.
The layers of $S_{f_x}$ coincide with layers of $S_f$ except
those are contained in $\beta_{uv}\bigcap\hat{A}_1$,
$\beta_{uv}\bigcap\hat{A}_3$ or $\beta_{uv}\bigcap\hat{A}_2 +s$. Therefore,
any representative $\pi$ lies in $S_{\hat{A}_1}\times
S_{\hat{A}_3}$ and has the form $\pi=\pi_1\pi_3$.
Consider two summands $\rho_1(\pi)\sigma\tau_1 x^{-1}\rho_2(\pi)^{-1}$ 
and 
$\rho_1(\pi')\sigma\tau_2 y^{-1}\rho_2(\pi')^{-1}$. If they are same we
have $\sigma^{-1}\rho_1(a)\sigma=\tau_2 y^{-1}\rho_2(a)^{-1}x\tau_1^{-1}$,
where $a=a_1a_3=\pi'^{-1}\pi=\pi_1'^{-1}\pi_1\pi_3'^{-1}\pi_3$.  
As in Lemma 3 from \cite{zub1} we see that
$\sigma^{-1}\rho_1(a)\sigma=\sigma^{-1}a_1\sigma\in H$ since there are no
any passive layers. In particular, $a\in 
\rho_1^{-1}(H^{\sigma})\bigcap 
\rho_2^{-1}(xHy^{-1})$. If $x=y$ then $a=
\pi'^{-1}\pi\in S_{f_x}$, that is the summands of $S(x)$ do not appear twice. 
The case $x\neq y$ means that
we have two equal summands from $S(x)$ and $S(y)$. As above,
it follows that $\rho_1(a)=\sigma h_1\sigma^{-1}$ and $\rho_2(a)=xh_2y^{-1}, 
h_1,h_2\in H$, $a\in S_f$.
For any $\pi\in S_f, \tau\in H$ we obtain 

$$\rho_1(\pi)\sigma\tau
x^{-1}\rho_2(\pi)^{-1}=\rho_1(\pi)\rho_1(a)\sigma h_1^{-1}\sigma^{-1}\sigma\tau
x^{-1}x h_2 y^{-1}\rho_2(a)^{-1}\rho_2(\pi)^{-1}=$$
$$=\rho_1(\pi a)\sigma h_1^{-1}\tau
h_2 y^{-1}\rho_2(\pi a)^{-1}.$$
       
\noindent In other words, $S(x)$ and $S(y)$ consist of the same summands. 
This concludes the proof for ordinary $A_1$-prelayers. 

$({\bf ii})$ From now on one can assume that any ordinary layer 
$\beta_{uv}$ satisfies
either $\sigma(\beta_{uv})\subseteq\hat{A}_1$ or $\sigma(\beta_{uv})\subseteq
\hat{A}_2\bigcup\hat{A}_3$.    

Consider the case $\sigma(\beta_{uv})\subseteq
\hat{A}_2\bigcup\hat{A}_3$. Denote by $\gamma$ the set $\beta_{uv}\bigcap     
\hat{A}_1$. We work with sums $S'(x)$. It is clear that all $A_1$ and 
$A_3$-layers of $S'_{f_x}$ and $S_f$ are the same. In particular, any 
representative $\pi\in S_f/S'_{f_x}$ can be choosen in $S_{\hat{A}_2}$
and $\rho_2(\pi)=1\in H$. An equation 
$\rho_1(\pi)\sigma x\tau_1=\rho_1(\pi')\sigma y\tau_2$
takes place iff $a=\pi'^{-1}\pi\in S_f\bigcap S_{\hat{A}_2}$ satisfies
$\rho_1(a)\in \sigma yH(\sigma x)^{-1}$. It remains to repeat the final
computations from $({\bf i})$. 
The same arguments work in the case
when $\gamma$ is any $A_3$-prelayer of $\beta_{uv}$. Therefore,
one can assume that either $\beta_{uv}\subseteq\hat{A}_1$ or
$\beta_{uv}$ coincides with its $A_3$-prelayer.  

$({\bf iii})$ We consider $\gamma$ which is a 
$A_2$-prelayer contained in 
$\beta_{uv}$ and work with sums $S(x)$. 
For any $x\in S_{\bar{\beta}}$
the layers of $S_{f_x}$ coincide with the layers of $S_f$ except
those which are contained in
$\beta_{uv}+s$ (if $\beta_{uv}\subseteq\hat{A}_2$) or in 
$\beta_{uv}$ (if $\beta_{uv}\subseteq\hat{A}_3$). In
particular, any representative $\pi\in S_f/S_{f_x}$ can be choosen in
$S_{\hat{A}_3}$ and $\rho_1(\pi)=1\in H^{\sigma}$. 
An equation 
$\sigma\tau_1 x^{-1}\rho_2(\pi)^{-1}=\sigma\tau_2 y^{-1}\rho_2(\pi')^{-1}$
holds iff the element $a=\pi^{-1}\pi'\in S_f\bigcap S_{\hat{A}_3}$ 
satisfies $\rho_2(a)\in xHy^{-1}$. It remains to refer to the  
final computations from $({\bf i})$ again.

$({\bf iv})$ Now, we consider the case $\beta_{uv}\subseteq\hat{A}_1, 
\sigma(\beta_{uv})\subseteq\hat{A}_2\bigcup\hat{A}_3$. We extract 
a $A_2$-prelayer $\gamma\subseteq\beta_{uv}$ and work with sums $S(x)$. 
It is clear that
only some $A_1$-layers of $S_{f_x}$ are different from $A_1$-layers of $S_f$.
Moreover, all of them are sublayers of $A_1$-layers of $S_f$ contained in
$\beta_{uv}$. Thus all representatives $\pi\in S_f/S_{f_x}$
can be choosen in $S_{\beta_{uv}}$. In particular, $\rho_i(\pi)=\pi, i=1,2$.
It is obvious that this case is the same as $({\bf i})$.

$({\bf v})$ Finally, let $\sigma(\beta_{uv})\subseteq\beta_{fg}
\bigcap\hat{A}_1$. 
Let $\gamma$ is a $A_1$-layer of $S_f$ or a $A_3$-prelayer 
belonging to $\beta_{uv}$. We work 
with  sums $S'(x)$. Only some $A_1$-layers of $S'_{f_x}$ are different 
from $A_1$-layers of $S_f$ and all of them are sublayers of $\sigma(\beta_{uv})
\bigcap\beta_{fg}$. Thus all representatives $\pi\in S_f/S'_{f_x}$ can be 
choosen 
in $S_{\beta_{fg}}\bigcap S_{\hat{A}_1}\leq H$. As above, $\rho_i(\pi)=\pi,
i=1,2$. An equation $\pi\sigma x\tau_1\pi^{-1}=\pi'\sigma y\tau_2\pi'^{-1}$
holds iff $(\sigma y)^{-1} a\sigma x=\tau_2 a\tau_1^{-1}
\in H$, where $a=\pi'^{-1}\pi$. The final computations are already obvious. 
The proposition is proved.

Now everything is prepared to prove the main theorem. 
Let $\beta_{uv}$ be an ordinary layer and $\sigma(\beta_{uv})\subseteq
\hat{A}_2\bigcup\hat{A}_3$. More precisely, suppose that 
$\sigma(\beta_{uv})\subseteq (\sigma(\beta_{fg})-s)\bigcap\hat{A}_2$. 
The case $\sigma(\beta_{uv})+s \subseteq
\sigma(\beta_{fg})\bigcap\hat{A}_3$ can be checked in the same way. 
It is possible that $\beta_{fg}=\beta_0$.  
Denote by ${\bf X}$ a variable corresponding to the $A_2$-layer
$\sigma(\beta_{uv})\bigcap(\sigma(\beta_{fg}) 
-s)$.
Using the contracting rules we see that the right hand side neighbor of 
$\overline{{\bf X}}$
in the records of all summands $tr^*(\sigma\tau)$ is either
${\bf Y}$ or $\overline{{\bf Y}}$, where the variable ${\bf Y}$ corresponds 
to the 
$A_1$-layer $\beta_{uv}$ or to the $A_3$-layer $\beta_{uv}+s \ 
(\beta_{uv})$. 
More precisely, when $\beta_{uv}\subseteq\hat{A}_3$ we have
a product $p=\overline{{\bf X}} \ \overline{{\bf Y}}$ but in other cases -- 
$p=\overline{{\bf X}}{\bf Y}$. Notice that the path $p$ is closed iff 
$\beta_{uv}=\beta_{fg}\pm s$. 

One can represent $z$ as $z=\sum_{I, K}\alpha_{I, K}\sigma_{i_1}
(p_1)^{k_1}\ldots\sigma_{i_l}(p_l)^{k_l}$, where $I=\{i_1,\ldots , i_l
\}, \\
K=\{k_1,\ldots , k_l\}$ are collections of indices, 
$p_1, \ldots, p_l$ are (not necessary different) primitive cycles
(if some $p_k, p_m$ are equal then we suppose that $i_k\neq i_m$) and
for any $I, K$ we have $i_1 k_1\mid p_1\mid +\ldots + i_l k_l\mid p_l\mid
=r$. Up to possible repetitions among $p_1,\ldots , p_l,$ 
the set $\{p_1,\ldots , p_l\}$
coincides with the set of all primitive cycles belonging to at least one
summand $tr^*(\sigma\tau)$ of $z$ (see \cite{zub7, don2}). Since $p$ is not a
proper power, even if it is a cycle, we see that each $p_i$ either 
contains $p$ (it is possible that $p_i$ contains $p$ more than one time) or 
does not contain any arrow belonging to $p$.

We list all possibilities for origins or ends of $p$ as follows.

\begin{enumerate}

\item If $\beta_{uv}\subseteq\hat{A}_1$ then $t(p)=a^*, j_q, \ i(p)=b, i_{q'}; 
a, b\in V_{ord}, q, q'\in\Omega$.

\item If $\beta_{uv}\subseteq\hat{A}_2$ then $t(p)=a^*, j_q, \ i(p)=j_{q'}; 
a\in V_{ord}, q, q'\in\Omega$.

\item If $\beta_{uv}\subseteq\hat{A}_3$ then $t(p)=a^*, j_q, i(p)=b^*, j_{q'};
a, b\in V_{ord}, q, q'\in\Omega$.

\end{enumerate}

In the case $i(p)=b, t(p)=a^*, a, b\in V_{ord}$ we construct a new quiver 
$\tilde{Q}$ with the vertex set $\tilde{V}=V\bigcup\{a^*\}$ and the arrow set
$\tilde{A}=A\bigcup\{\tilde{p}\}$, where $\tilde{p}$ is a new arrow having
the same origin and end as the path $p$. In other words, the difference
between $Q$ and $\tilde{Q}$ is that the set of couples $\{i_q, j_q\}$ is 
completed by the new couple $\{a, a^*\}$ but $\tilde{V}_{ord}=V_{ord}\setminus
\{a\}$. For the sake of convenience we
introduce a symbol $q_0$ such that $a=i_{q_0}, a^*=j_{q_0}$ and
$\tilde{\Omega}=\Omega\bigcup\{q_0\}$. In all other cases it is not
necessary to add new vertices but only new arrows. For example, if $i(p)=j_q,
t(p)=a^*, a\in V_{ord}, q\in\Omega$, then $\tilde{V}=V, \tilde{A}=A\bigcup
\{\tilde{p}\}$, where $i(\tilde{p})=a, t(\tilde{p})=i_q$. Notice that $
\tilde{V}^{(d)}=V^{(d)}$ but $\tilde{A}^{(d)}$ is different from $A^{(d)}$
whenever $i(p)\in V_{ord}, t(p)\in V_{ord}^*$.

If $\tilde{V}=V$ we leave the same dimension vector ${\bf t}$ but
if $\tilde{V}\neq V$ then we replace it by $\tilde{{\bf t}}=(\ldots , d_a,
d_a^*,\ldots)$. It is clear that the representation space $R(\tilde{Q},
{\bf t}) \ (R(\tilde{Q}, \tilde{{\bf t}}))$ contains the space $R(Q, {\bf t})$
as a direct summand. Thus $J(Q, {\bf t})\subseteq J(\tilde{Q}, {\bf t})$ 
\ ($J(Q, {\bf t})\subseteq J(\tilde{Q}, \tilde{{\bf t}}))$ and $J(Q)\subseteq
J(\tilde{Q})$.

Replace all occurences of $p$ in $z$ by $\tilde{p}$. We get some
$\tilde{z}\in J(\tilde{Q})$. Since matrices $X_{{\bf d}}, Y_{{\bf d}}$
appear only in the product $p$, it is easy to prove that $z\in T(Q, {\bf t})$ 
iff $\tilde{z}\in T(\tilde{Q}, {\bf t}) \ (\tilde{z}\in T(\tilde{Q}, 
\tilde{{\bf t}}))$. Using the induction hypothesis we obtain

$$\tilde{z}=\sum_{i\in V_{ord}, u\geq d_i+1} f_{i, u}\sigma_u(h_{i, u})+
\sum_{i\in V_{ord}, v\geq d_i+1, 2s\leq v} f_{i, v, s}\sigma_{v, s}(
h^{(1)}_{i, v, s}, h^{(2)}_{i, v, s}, h^{(3)}_{i, v, s})+$$  
$$+\sum_{q\in\Omega, v\geq d_q+1, 2s\leq v} f_{q, v, s}\sigma_{v, s}(
h^{(1)}_{q, v, s}, h^{(2)}_{q, v, s}, h^{(3)}_{q, v, s}),$$

or

$$\tilde{z}=\sum_{i\in \tilde{V}_{ord}, u\geq d_i+1} f_{i, u}
\sigma_u(h_{i, u})+
\sum_{i\in \tilde{V}_{ord}, v\geq d_i+1, 2s\leq v} f_{i, v, s}
\sigma_{v, s}(h^{(1)}_{i, v, s}, h^{(2)}_{i, v, s}, h^{(3)}_{i, v, s})+$$
$$+\sum_{q\in\tilde{\Omega}, v\geq d_q+1, 2s\leq v} f_{q, v, s}\sigma_{v, s}(
h^{(1)}_{q, v, s}, h^{(2)}_{q, v, s}, h^{(3)}_{q, v, s}).$$
  
\noindent Here $f_{i, u}, f_{i, v, s}, f_{q, v, s}$ are some monomials from
$J(\tilde{Q})$. The elements
$h_{i, u}$ are incident to $i$ and $h^{(1)}_{i, v, s}, h^{(1)}_{q, v, s}$
are incident to $i$ and $i_q$ respectively. Furthemore, $h^{(k)}_{i, v, s},
h^{(k)}_{q, v, s}$ are passing from $i \ (i^*)$ and $i_q \ (j_q)$ to 
$i^* \ (i)$ and
$j_q \ (i_q)$ correspondingly, $k=2 \ (k=3)$. Replacing all occurences of 
$\tilde{p}$ by the product $p$ we complete the proof.

The case $\sigma(\beta_{uv})\subseteq\hat{A}_1\bigcap\beta_{fg}$ is the same.
As above, ${\bf X}$ is a variable corresponding to the $A_1$-layer 
$\sigma(\beta_{uv})\bigcap\beta_{fg}$. It is easy to see that the right hand 
side neighbor of ${\bf X}$ is ${\bf Y}$ or $\overline{{\bf Y}}$, where 
${\bf Y}$ corresponds to either
an $A_1$-layer $\beta_{uv}$ or to an $A_3$-layer $\beta_{uv}+s \ (\beta_{uv})$.
We leave to the reader to check all details.

It remains to consider the last case when there are no any ordinary layers.
It means that either $V=V_{ord}=\{1\}$ or $\Omega=\{q\}, V_{ord}=\{i_q\}$.
In both cases $S_{\bar{\beta}}=S_{\beta_0}=S_{{\cal T}}=S_{{\cal T}(x)}=
S_{{\cal H}}=S_{{\cal H}(x)}=S_r$, $x=1, q$. If $x=1$ then our quiver consists
of one loop incident to $1$ and the element $z$ is equal to 
$\sigma_r({\bf X}), r > d_1$. If $x=q$ then the quiver has two vertices 
$i_q, j_q$,
one loop incident to $i_q$, one arrow passing from
$i_q$ to $j_q$ and one arrow having opposite direction. It is clear that
$z=\sigma_{r, s}({\bf X}, {\bf Y}, {\bf Z})$. The theorem is proved.

\begin{rem}
It is possible to give a self-contained proof of this theorem which does not
use preliminary description of the free invariant algebra $J(Q)$ given in
\cite{zub7}. 
In fact,
one has to prove that every time when we replace given suitable 
generator $z$ by $\tilde{z}$, as above, we get a suitable generator 
again with respect 
to some other base group and quiver. But there are two reasons why I prefered 
the way used in this article. 

First, it is not obvious that the elements 
$\sigma_{r, s}( {\bf X}, {\bf Y}, {\bf Z})$ can be written as sums  
$\sum_{I, K}\alpha_{I, K}
\sigma_{i_1}(p_1)^{k_1}\ldots\sigma_{i_l}(p_l)^{k_l}$. Of course, referring to
\cite{zub7} we know that it is true, but how to get this expression for any
$\sigma_{r, s}({\bf X}, {\bf Y}, {\bf Z})$ directly? 

Second, it is not easy exercise to show that $\tilde{z}$ is a 
suitable generator 
and it requires a lot of case-by-case observations. For example, 
consider the case $\sigma(\beta_{uv})\subseteq(\sigma(\beta_0)-s)\bigcap
\hat{A}_2, \beta_{uv}=\sigma(\beta_{cd})\subseteq\hat{A}_1, t(p)=a^*,
i(p)=b$, where $\beta_{cd}$ is an ordinary layer and $a, b\in V_{ord}$.
The conditions of admissibility say that $\beta_{uv}\subseteq T(i_q)\bigcap
\hat{A}_1,
\sigma(\beta_{uv})\subseteq I(j_q)$, $\mu\subseteq {\cal T}(a)$ and
$\beta_{cd}\subseteq {\cal T}(b)$, where $\mu\subseteq\beta_0, \sigma(\mu)-s=
\sigma(\beta_{uv})$. One can check that

$$\tilde{z}=c(\sum_{\tau\in S_{\bar{\beta'}}}(-1)^{\tau}
\sigma'\tau, S_{f'}).$$

\noindent Here $S_{\bar{\beta'}}=S_{\bar{\beta}\setminus\beta_{uv}}, 
\sigma'\mid_{[1, r]
\setminus (\beta_{uv}\bigsqcup\beta_{cd}\bigsqcup\mu)}=\sigma$ but $\sigma'
\mid_{\mu}=\sigma-s, \sigma'\mid_{\beta_{cd}}=\sigma^2+s$. In other words, 
if $\bar{r'}$ is a multidegree of $\tilde{z}$ then $r'_{{\bf X}}=
r'_{{\bf Y}}=0, 
r'_{\tilde{p}}=r_{{\bf X}}=r_{{\bf Y}}=\mid\beta_{uv}\mid$, the variable 
$\tilde{p}$ 
corresponds to the $\tilde{A}_2$-layer $\sigma(\beta_{uv})=\sigma'(\mu)=
\sigma'(\beta_{cd})-s$. Moreover, $S'_0=S_{\hat{A}_1\setminus\beta_{uv}}\times
S_{\hat{A}_2}\times S_{\hat{A}_3}$ and $\rho'_1, \rho'_2 : S'_0\to
S_{[1, r]\setminus\beta_{uv}}$ are just restrictions of $\rho_1, \rho_2$.
It is clear that the layers of $S_{f'}=\rho_1^{'-1}(S_{\bar{\beta'}}
^{\sigma'})\bigcap\rho_2^{'-1}(S_{\bar{\beta'}})$ are layers of $S_f$ without
$A_1$-layer $\beta_{uv}$. One also has to check that
$\sigma'\in S_{[1, r]\setminus\beta_{uv}}$. We have $[1, r]\setminus\beta_{uv}
=\mu\bigsqcup\beta_{cd}\bigsqcup T$. Since $\sigma'$ acts injectively on the
subsets $\mu, \beta_{cd}, T$ it remains to prove that 
$\beta_{uv}$ does not intersect the set $\sigma'([1, r]\setminus\beta_{uv})$.
By definition

$$
\sigma'([1, r]\setminus\beta_{uv})=\sigma'(\beta_{cd})\bigcup\sigma'(\mu)
\bigcup\sigma(T)=(\sigma(\beta_{uv})+s)\bigcup\sigma(\beta_{uv})\bigcup
\sigma(T)=$$
$$=\sigma(\mu)\bigsqcup\sigma(\beta_{uv})\bigsqcup\sigma(T).$$

\noindent Now it is obvious because of $\beta_{uv}=\sigma(\beta_{cd})$.     
\end{rem}

\section{Applications}

The notation of mixed representations of quivers was introduced to generalize
Procesi-Razmyslov's theorem (briefly -- {\bf PRT})
for adjoint action invariants of orthogonal and
symplectic groups (see \cite{zub5, zub6, zub7}). 
In this section modulo the previous theorem we describe some approach
to this problem. In fact, the same  method works  
in much more general case of so-called {\it supermixed} representations of 
quivers (see for definitions \cite{zub7}). The procedure of computation of
generating invariants of supermixed representations of any quiver was
described in \cite{zub7}. To be precise, these invariants can be obtained by 
a specialization of invariants of mixed representations of another quiver 
(see Section 4 from \cite{zub7}). The next step is to get all defining
relations between them. We demonstrate how to do it in the principal case
of the diagonal actions of orthogonal and symplectic groups on several
matrices by conjugation. The general case can be reduced to the principal
one as in \cite{zub7}.

Consider the quiver $Q$ such that $V=\{1,2\}$ and $A=\{a_1,\ldots a_m,
b, c\}$ with $i(a_j)=t(a_j)=1, i(b)=t(c)=1, t(b)=i(c)=2, 1\leq j\leq m$.  
Let $V_{ord}=\emptyset, \Omega=\{q\}, i_q=1, j_q=2$. Any dimension vector
${\bf t}$ compatible with this partition has the form $(d, d^*)$.
It was proved in \cite{zub5} that for any $d$ we have a short exact sequence

$$0\to I_d\to J(Q, d)\to S_d\to 0.$$

\noindent Here $J(Q, d)=J(Q, (d, d^*)), S_d=K[M(d)^m]^{G_d}, G_d=O(d)$ or 
$G_d=Sp(d)$.
The ideal $I_d$ is equal to $T_d^{GL(d)}$, where $T_d$ is the ideal of
$K[R(Q, (d, d^*))]=K[R(Q, d)]$ generated by the coefficients of the matrices
$Y(b)Y(c)-E(d), Y(c)-\overline{Y(c)}$ \ (in the
symplectic case, the last matrix should be replaced by 
$Y(c)+\overline{Y(c)}$ if $char K\neq 2$, otherwise one has to fill 
its zero diagonal by the 
original coefficients $y_{kk}(c), k=1,\ldots, d$), $E(d)$ is an $d\times d$ 
unit matrix. 
The epimorphism $J(Q, d)\to S_d\to 0$ is induced by
the specialization $Y(b), Y(c)\mapsto E(d)$ \ (by $Y(b)\mapsto J(d), Y(c)
\mapsto -J(d)$ in the symplectic group case, where $J(d)$ is a matrix of the 
skew-symmetric bilinear form defining $Sp(d)$). Recall that in the
orthogonal group case we assume that $char K\neq 2$.        

First, we consider the orthogonal group case. Fix two integers 
$N, n, N > 
n$. For the dimension vectors ${\bf N}=(N, N^*), {\bf n}=(n, n^*)$ we
denote the non-standard (on the arrows $b$ and  $c$) specialization 
$p_{{\bf N}, {\bf n}}$ ($j_{{\bf N}, {\bf n}}$) just by $p_{N, n}$
(respectively -- by $j_{N, n}$). 
Denote by $E_d$ the subspace of
$K[R(Q, d)]$ which is generated by the polynomials 
$z_{ij}=\sum_{1\leq k\leq d} 
y_{ik}(b)y_{kj}(c)-\delta_{ij}, u_{kl}=y_{kl}(c)-y_{lk}(c), i,j, k, l=1,
\ldots , d, k < l$. It is easy to see that 
$Z^g=g^{-1}Zg, U^g=g^{-1}U\bar{g}^{-1}$, 
where $Z=(z_{ij}), U=\frac{1}{2}(Y(c)-\overline{Y(c)}), g\in GL(d)$. 

Since the elements $z_{ij}, u_{kl}$ form a regular sequence in $K[R(Q, d)]$, 
we have the Koszul resolution (see Lemma1.3(b) \cite{don7} or \cite{mats})

$$\ldots\to\Lambda^k(E_d)\otimes K[R(Q, d)]\to\ldots\to
E_d\otimes K[R(Q, d)]\to T_d\to 0.$$

For any $k$ denote by $\Delta_{d, k}$ the image of $\Lambda^{k+1}(E_d)
\otimes K[R(Q, d)]$ in $\Lambda^k(E_d)\otimes K[R(Q, d)]$. Notice that 
$g.f.d.(S^l(E_d))=0$ for any $l\geq 0$ \cite{kur1, kur2}. 
Since $g.f.d(
\Lambda^k(E_d)\leq k-1$ (\cite{don7}, Corollary
1.2(d)) we obtain $g.f.d(\Lambda^k(E_d)\otimes 
K[R(Q, d)])\leq k-1$ because of
$g.f.d.(K[R(Q, d)])=0$ \cite{zub7} and
$g.f.d(\Delta_{d, k})\leq k, k=1,2,\ldots,$ by Lemma 1.1.

\begin{pr} The specialization $p_{N, n}$ induces an epimorphism $I_N\to
I_n$.
\end{pr}

Proof. It is easy to see that 
$p_{N, n}
(T_N)=T_n, p_{N, n}(\Lambda^k(E_N)\otimes K[R(Q, N)])=
\Lambda^k(E_n)\otimes K[R(Q, n)]$ and $p_{N, n}(I_N)
\subseteq I_n$. 
Moreover, for any $k\geq 1$ we have
$p_{N, n}(\Delta_{N, k})\subseteq\Delta_{n, k}$. Applying Lemma 1.2
for the homomorphism of pairs $(p_{N, n}, j_{N, n})$ one can extract the 
following
commutative fragments of the corresponding long diagrams (we omit the first
arguments $GL(N), GL(n)$ from all cohomology groups because of they can be 
easily recovered by referring to subindices)

$$\begin{array}{ccccccc}
(E_N\otimes K[R(Q, N)])^{GL(N)} & \to & I_N & \to &
H^1(\Delta_{N, 1})& \to & 0 \\
\downarrow & & \downarrow & & \downarrow & &  \\
(E_n\otimes K[R(Q, n)])^{GL(n)} & \to & I_n &
\to & H^1(\Delta_{n, 1})& \to & 0 ,
\end{array}
$$

\noindent and (for $k=2,3, \ldots$)

$$\begin{array}{ccccccc}
H^{k-1}(\Lambda^k(E_N)\otimes K[R(Q, N)])& \to &
H^{k-1}(\Delta_{N, k-1}) & \to &
H^k(\Delta_{N, k})& \to & 0 \\
\downarrow & & \downarrow & & \downarrow & &  \\
H^{k-1}(\Lambda^k(E_n)\otimes K[R(Q, n)])& \to &
H^{k-1}(\Delta_{n, k-1}) & \to &
H^k(\Delta_{n, k})& \to & 0 .
\end{array}
$$ 

Assume  that all $H^{k-1}(\Lambda^k(E_N)\otimes K[R(Q, N)])
\to H^{k-1}(\Lambda^k(E_n)\otimes K[R(Q, n)])$ are epimorphisms, $k
\geq 1$.
Then $I_N\to I_n$ is an epimorphism iff $H^1(\Delta_{N, 1})\to
H^1(\Delta_{n, 1})$ is. Regarding to the next diagram we see that
$H^1(\Delta_{N, 1})\to
H^1(\Delta_{n, 1})$ is an epimorphism iff $H^2(\Delta_{N, 2})
\to H^2(\Delta_{n, 2})$ is and so on. But $\Delta_{n, k}=0$ for 
sufficiently large $k$. Therefore, all we need is to prove that 
$H^{k-1}(\Lambda^k(E_N)\otimes K[R(Q, N)])\to H^{k-1}(\Lambda^k(E_n)
\otimes K[R(Q, n)])$ is an epimorphism for any $k\geq 1$. 

Consider some collection of symmetric powers 
$S^{l_1}(E_d), \ldots, S^{l_s}(E_d), 
d=N, n$. 
For the sake of simplicity denote
$K[R(Q, d)]$ by $A_d$ and $S^{l_1}(E_d)\otimes\ldots\otimes S^{l_s}(E_d)$ by
$B_d$. By the definition of $E_d$ we have 

$$B_d=
\oplus_{0\leq i_1\leq l_1,\ldots , 0\leq i_s\leq l_s} S^{i_1}(Z_d)\otimes
S^{l_1-i_1}(U_d)\otimes\ldots\otimes S^{i_s}(Z_d)\otimes S^{l_s-i_s}(U_d).$$

\noindent Here $S^a(Z_d)=K[Z_d](a), S^b(U_d)=K[U_d](b)$ are homogeneous 
components of 
the polynomial algebras generated by the coefficients of $Z_d$ and $U_d$ 
correspondingly.

The following lemma completes the proof of Proposition 3.1.

\begin{lm}
The homomorphism 
$$H^{k-1}(\Lambda^k(E_N)\otimes A_N\otimes B_N)\to 
H^{k-1}(\Lambda^k(E_n)\otimes A_n\otimes B_n)$$ 
induced by $p_{N, n}$ is an epimorphism for any $k\geq 1$. 
\end{lm}

Proof. The space
$(E_N\otimes A_N\otimes B_N)^{GL(N)}$ can be regarded as a sum of
homogeneous  components
of the invariant algebra of supermixed representations of the other quiver 
$Q'$ with the same vertex set as $Q$ but with $s+1$ new loops incident to
the vertex $1$ and $s+1$ new arrows passing from the vertex $2$ to the vertex 
$1$. The spaces of the linear maps which belong to the last arrows are the 
spaces of
skew-symmetric matrices corresponding to $s+1$ copies of $U$. The spaces 
belonging to  the new loops are the spaces of square matrices
corresponding to $s+1$ copies of $Z$.
For more detailed explanations we refer to  
\cite{zub7}. Finally, the restriction $p_{N, n}$ on 
$(E_N\otimes A_N\otimes B_N)^{GL(N)}$ can be
identified with the standard specialization of the corresponding invariant 
algebras 
which is a homogeneous epimorphism (see \cite{zub7}, Section 4). 
Thus yields the case $k=1$.

Let $k > 1$. We have an exact sequence (it is also a partial case of the 
Koszul resolution \cite{don7, mats, bur})

$$0\to\Lambda^k(E_d)\to \Lambda^{k-1}(E_d)\otimes S^1(E_d)\to
\ldots\to\Lambda^1(E_d)\otimes S^{k-1}(E_d)\to S^k(E_d)\to 0.$$

\noindent Tensoring by $C_d=A_d\otimes B_d$ we get an exact sequence

$$0\to\Lambda^k(E_d)\otimes C_d\to 
\Lambda^{k-1}(E_d)\otimes S^1(E_d)\otimes C_d\to$$
$$\ldots\to\Lambda^1(E_d)\otimes S^{k-1}(E_d)\otimes C_d
\to S^k(E_d)\otimes C_d\to 0.$$

As above, denote by $\Delta_{d, i}$ the image of $\Lambda^{i+1}(E_d)\otimes
S^{k-i-1}(E_d)\otimes 
C_d$ in $\Lambda^i(E_d)\otimes S^{k-i}(E_d)\otimes C_d$. 
Repeating the previous arguments we obtain a collection of commutative 
diagrams with exact top and bottom rows

$$\begin{array}{ccccccc}
H^{k-2}(\Lambda^{k-1}(E_N)\otimes S^1(E_N)\otimes C_N) & 
\to & H^{k-2}(\Delta_{N, k-2}) & \to &
H^{k-1}(\Lambda^k(E_N)\otimes C_N) & 
\to & 0 \\
\downarrow & & \downarrow & & \downarrow & &  \\
H^{k-2}(\Lambda^{k-1}(E_n)\otimes S^1(E_n)\otimes C_n) & 
\to & H^{k-2}(\Delta_{n, k-2}) & \to &
H^{k-1}(\Lambda^k(E_n)\otimes C_n) & 
\to & 0,
\end{array}
$$

\noindent and (for $i=3, \ldots$)

$$\begin{array}{ccccccc}
H^{k-i}(\Lambda^{k-i+1}(E_N)\otimes S^{i-1}(E_N)\otimes C_N) & 
\to & H^{k-i}(\Delta_{N, k-i}) & \to &
H^{k-i+1}(\Delta_{N, k-i+1}) & 
\to & 0 \\
\downarrow & & \downarrow & & \downarrow & &  \\
H^{k-i}(\Lambda^{k-i+1}(E_n)\otimes S^{i-1}(E_n)\otimes C_n) & 
\to & H^{k-i}(\Delta_{n, k-i}) & \to &
H^{k-i+1}(\Delta_{n, k-i+1}) & 
\to & 0.
\end{array}
$$

By induction hypothesis $H^{k-2}(\Lambda^{k-1}(E_N)\otimes 
S^1(E_N)\otimes C_N)\to H^{k-2}(\Lambda^{k-1}(E_n)\otimes 
S^1(E_n)\otimes C_n)$ is an epimorphism. Thus $H^{k-1}(\Lambda^k(E_N)\otimes 
C_N)\to H^{k-1}(\Lambda^k(E_n)\otimes C_n)$ is an
epimorphism iff $H^{k-2}(\Delta_{N, k-2})\to H^{k-2}(\Delta_{n, k-2})$ is.
It is clear that the next typical step is to show that 
$H^{k-i+1}(\Delta_{N, k-i+1})
\to H^{k-i+1}(\Delta_{n, k-i+1})$ is an epimorphism iff
$H^{k-i}(\Delta_{N, k-i})\to H^{k-i}(\Delta_{n, k-i})$ is. But for $i=k$
it is obviously the base of induction. The lemma and the proposition are 
proved.

Denote by $\tilde{p}_{N, n}$ the standard specialization $J(Q,N)
\to J(Q, n)$. Using Amitsur's formulae \cite{am} we see that $p_{N, n}$
takes any $\sigma_j(m)$ to either $\sigma_j(\tilde{p}_{N, n}(m))$ or to 
$\sum_{0\leq l\leq j}\sigma_l(\tilde{p}_{N, n}(m))\sigma_{j-l}(E(N, n))$, 
where
$m\in Q^{(d)}$ (as above we identify monomials with pathes in $Q^{(d)}$) and
$E(N, n)$ is a $N\times N$ matrix with $1$-s on the diagonal except the 
first $n$ places and $0$-s on all other places. In particular, we get

\begin{lm}
For any $f\in J(Q, N)(r)$ $p_{N, n}(f)=\tilde{p}_{N, n}(f) + (\mbox{summands of
degree} < r)$.
\end{lm} 
 
If $A$ is a graded algebra we denote by $A^{(r)}$ the sum 
$\oplus_{0\leq i\leq r}A(i)$. So the algebra $A$ turns to a filtred algebra 
with filtration $A^{(0)}\subseteq A^{(1)}\subseteq\ldots\subseteq A^{(r)}
\subseteq\ldots$. For any ideal $I$ of $A$ denote $I\bigcap A^{(r)}$ 
by $I^{(r)}, r=0,1,2,\ldots$.

\begin{lm}
For any $r\geq 0$ we have $p_{N, n}(J(Q, N)^{(r)})=J(Q, n)^{(r)}$. 
Moreover, if $N'\geq N\geq r$ then $p_{N', N}\mid_{J(Q, N')^{(r)}}$ 
is an isomorphism.   
\end{lm}

Proof. Use induction on $r$ and Theorem 
2.1 from \cite{zub7}. If $f\in J(Q, N')^{(l)}\setminus J(Q, N')^{(l-1)}, l\leq
r,$ then by Lemma 3.2 $f=f_l + (\mbox{summands of degree} <l)$, where 
$f_l$ is the non-zero 
$l$-th homogeneous component of $f$. Thus $p_{N', N}(f)=\tilde{p}_{N', N}
(f_l) + (\mbox{summands of degree} <l)$ and $\tilde{p}_{N', N}(f_l)\neq 0$
\cite{zub7} (see the remark after Theorem 1 or Proposition 2.2 from
\cite{zub7}). The lemma is proved.

\begin{rem}
By \cite{zub5} the algebra $S_d$ is generated by the elements $\sigma_j(p),
1\leq j\leq d$, where $p$ is an arbitrary product of matrices 
$Y(a_i), \overline{Y(a_i)}, 1\leq i\leq m$. 
Thus it obviously follows that $p_{N, n}(S_N(r))=S_n(r)$
for any $N > n, r\geq 0$. Moreover, it is not hard to prove that
the space $S_{n}^{(r)}$ is covered by $J(Q, n)^{(2r)}$. Summarizing all
previous statements one can say that for any $N\geq n$ in the 
commutative diagram

$$\begin{array}{ccccccccc}
0 & \to & I_N &\to & J(Q, N) &\to & S_N & \to & 0 \\
  &         &\downarrow &  & \downarrow &        & \downarrow & & \\
0 & \to & I_n &\to & J(Q, n) &\to & S_n & \to & 0  
\end{array}
$$

\noindent all vertical arrows are epimorphisms.
\end{rem}

\begin{lm}
For a fixed $r$ and $N'\geq N\geq r,$ 
$p_{N', N}\mid_{S_{N'}(r)}$ is an isomorphism.
\end{lm}
 
Proof. It is equivalent to prove that $\dim S_{N'}(r)=\dim S_N(r)$. Since
$K[M(n)^m](r)$ is an $O(n)$-module with good filtration for any $n, r$
\cite{zub3} we see that
$\dim S_n(r)=\dim K[M(n)^m](r)^{O(n)}$ equals the multiplicity of a trivial
module and does not depend on the characteristic of the ground field $K$.
In fact, the formal character of $K[M(n)^m](r)$ 
as well as its representation as a sum of formal characters of induced
modules does not depend on the characteristic of $K$
(see \cite{jan, don2, don3, don4, don5, don6, zub1, zub2, zub3}
for more explanations).
In particular, one can suppose that $char K=0$. If $S_{N'}(r)\to
S_N(r)\to 0$ is not an isomorphism that there is $f\in S_{N'}(r)\setminus 
0$ such that $p_{N', N}(f)=0$. It remains to replace $f$ by its complete
linearization and refer to \cite{pr}. The lemma is proved.

We have the countable set of spectrums $\{ J(Q, n)^{(r)}, 
p_{N, n}\mid N\geq n\}, r=0,1,2,\ldots$. Denote by $J'(Q)^{(r)}$ 
the inverse limit of $r$-th
spectrum. By Lemma 3.3 one can identify $J'(Q)^{(r)}$ 
with $J(Q, N)^{(r)}$ for sufficiently large $N$. In particular,
we have an inclusion $J'(Q)^{(r)} \to J'(Q)^{(r')}$ for any $r'\geq r$. 
Denote by $J'(Q)$ the direct limit of the spectrum consisting of all 
$J'(Q)^{(r)}$ and inclusions defined above. 

\begin{lm}
The graded algebra $\gr J'(Q)=\oplus_{r\geq 0}J'(Q)^{(r+1)}/J'(Q)^{(r)}$ is
isomorphic to $J(Q)$. 
\end{lm}  

Proof. It is sufficient to notice that the following diagram

$$\begin{array}{ccccc}
J(Q, N')^{(r+1)}/J(Q, N')^{(r)} & 
\to & J(Q, N')^{(r+1)}/J(Q, N')^{(r)} & = &
J(Q, N')(r+1) \\
p_{N', N}\downarrow & & \tilde{p}_{N', N}\downarrow & 
& \tilde{p}_{N', N}\downarrow \\
J(Q, N)^{(r+1)}/J(Q, N)^{(r)} & 
\to & J(Q, N)^{(r+1)}/J(Q, N)^{(r)} & = &
J(Q, N)(r+1)
\end{array}
$$

\noindent is commutative for any $N'\geq N\geq r+1$. 
Here the first two vertical maps are induced by $p_{N', N}$
and $\tilde{p}_{N', N}$ respectively. The horizontal maps are natural 
identifications. We leave checking of all rest details to the reader.

\begin{rem}
As it was noticed in the introduction one can prove that $J(Q)$ is a polynomial
algebra with homogeneous free generators. 
They can be choosed as $\sigma_j(p)$, where
$p$ runs over all primitive cycles.  In particular, $J'(Q)\cong J(Q)$!
\end{rem}

By Lemma 3.3 we have an epimorphism $J'(Q)\to J(Q, n)$. Denote by $T'(Q, n)$
the kernel of this epimorphism. Since $p_{N', N}$ coincides with 
$\tilde{p}_{N', N}$
on $S_{N'}$ the definitions of a free  
algebra of orthogonal invariants as a filtred or graded algebra are 
the same. We denote this algebra by
$S$. As above any homogeneous component $S(r)$ can be naturally identified 
with $S_N(r)$ for sufficiently large $N$ and we have an epimorphism
$S\to S_n$ with a kernel $K_n$.     

Now our aim is to describe the generators of $K_n$ as we declared at the
beginning of this section.
Let $f\in K_n^{(r)}$. Without loss of generality one can
assume that $f\in S_N^{(r)}, N >> r$. 

\begin{lm}
For sufficiently large $N, r', N\geq r'\geq r,$ there is $t\in T'(Q, n)^{(r')}=
T'(Q, n) \bigcap J(Q, N)^{(r')}$ 
such that the epimorphism
$J(Q, N)\to S_N$ takes $t$ to $f$.
\end{lm}

Proof. As we noticed in Remark 3.1 there is $f'\in J(Q, N)^{(2r)}$ such that
$f$ is the image of $f'$. Thus the image of $f''=p_{N, n}(f')$ in $S_n$ equals
zero, that is $f''\in I_n^{(2r)}$. By Proposition 3.1 
for sufficiently large $r'\geq 2r$
there is $f'''\in I_N^{(r')}$ such that $p_{N, n}(f''')=f''$.
In particular, $p_{N, n}(f'-f''')=0$. Increasing $N$ one can assume that 
$N\geq r'$. It remains to take a preimage of $f'-f'''$, say 
$t$, in $J(Q, N)^{(r')}$. The lemma is proved.

Let $f\in T(Q, n)^{(r)}$, that is $f\in
J(Q, N)^{(r)}, N >> r$. By Lemma 3.2 $p_{N, n}(f)=g'\in J(Q, n)^{(r-1)}$ and by
Lemma 3.3 one can choose an element $g\in J(Q, N)^{(r-1)}$ such that 
$p_{N, n}(g)=g'$. Let $g_{r-1}$ be a $(r-1)$-th homogeneous
component of $g$. Again, by Lemma 3.2 $p_{N, n}(f-g_{r-1})=h'\in
J(Q, n)^{(r-2)}$ and one can repeat the previous step. 
After $r$ steps like above we obtain some
$\tilde{f}\in T'(Q, n)^{(r)}$ such that it has $r$-th homogeneous component
$\tilde{f}_r$ coincided with $f_r$ -- $r$-th homogeneous component of $f$.  
Notice that if $f$ does not depend on $Y(b), Y(c)$ then $\tilde{f}=f$. 
        
\begin{lm}
The ideal $T'(Q, n)$ is generated by the elements $\tilde{f}$, where $f$ runs
over the set of generators of $T(Q, n)$ from Theorem 2.
\end{lm}

Proof. Fix some $N >> r$ and consider an element $f\in T'(Q, n)\bigcap 
J(Q, N)^{(r)}$. If $f_r$ is the $r$-th homogeneous component of $f$ then by
Lemma 3.2 we get
$p_{N, n}(f)=\tilde{p}_{N, n}(f_r) + (\mbox{summands of degree} < r)=0$. 
In particular, $\tilde{p}_{N, n}(f_r)=0$, that is $f_r\in T(Q, n)(r)$ and
$f_r$ can be represented as $\sum h_ig_i$, where any $h_i$ is
a homogeneous element from $J(Q, N)$ and $g_i$ is a homogeneous 
component of some generator from Theorem 2.
It is clear that an element 
$t=\sum h_i\tilde{g}_i$ lies in $T'(Q, n)$ and has the same $r$-th homogeneous
component as $f$, that is $f-t\in T'(Q, n)\bigcap 
J(Q, N)^{(r-1)}$. Induction on $r$ completes the proof.

The symplectic group case can be treated in the same way up to some change in 
the initial notations. To be precise, in this case $N=2M, n=2m$ and $p_{N, n}$
must be redefined as 

$$p_{N, n}(y_{ks}(b))=\left\{\begin{array}{l}
y_{ks}(b), \ \mbox{if} \ M-m+1\leq k, s\leq M+m, \\
1, \ \mbox{if} \ k+s=N+1, 1\leq k\leq M-m, \\
-1, \ \mbox{if} \ k+s=N+1, M+m+1\leq k\leq N, \\
0, \ \mbox{otherwise},
\end{array} \right.
$$

$$p_{N, n}(y_{ks}(c))=\left\{\begin{array}{l}
y_{ks}(c), \ \mbox{if} \ M-m+1\leq k, s\leq M+m, \\
-1, \ \mbox{if} \ k+s=N+1, 1\leq k\leq M-m, \\
1, \ \mbox{if} \ k+s=N+1, M+m+1\leq k\leq N, \\
0, \ \mbox{otherwise}.
\end{array} \right.
$$

\noindent On the rest variables $p_{N, n}$ acts by the 
{\it symplectically standard} rule

$$p_{N, n}(y_{ks}(a_i))=\left\{\begin{array}{l}
y_{ks}(a_i), \ \mbox{if} \ M-m+1\leq k, s\leq M+m, \\
0, \ \mbox{otherwise},
\end{array} \right.
\ , 1\leq i\leq m.$$   

Similarly, $\tilde{p}_{N, n}$  is symplectically standard on all variables.
The homomorphism $j_{N, n}$ is redefined as

$$(j_{N, n}(g))_{ks}=\left\{\begin{array}{l}
g_{k-M+m, s-M+m}(b), \ \mbox{if} \ M-m+1\leq k, s\leq M+m, \\
1, \ \mbox{if} \ k=s, 1\leq k\leq M-m \ \mbox{or} \ M-m+1\leq k\leq N \\
0, \ \mbox{otherwise},
\end{array} \right. \  , g\in GL(n). 
$$   

Denote by the same symbol $i_{N, n}$ as above
the morphism dual to $p_{N, n}$.
A free invariant algebra of symplectic invariants as well as a kernel
of an epimorphism of this algebra onto the algebra of symplectic invariants of
$m$ $n\times n$ matrices will be denoted by the same symbols $S, K_n$. 
It is easy to see that the invariant algebra $J(Q)$ remains the same even if we
replace standard specializations by symplectically standard. The proof of
Proposition 3.1 and all consequent lemmas can be word by word repeated.
Notice that the correspondence $f\mapsto \tilde{f}$ from Lemma 3.7 is 
different from the orthogonal group case because of one has to replace the 
matrix $E(N, n)$ by $J(N, n)=i_{N, n}(0)$.
Summarizing we have

\begin{pr}
The ideal $K_n$ (both orthogonal or symplectic) is generated by the images of 
the elements $\tilde{f}$ from
Lemma 3.7.
\end{pr}

\section{Concluding remarks}

Proposition 3.2 gives only some procedure to compute the generators of $K_n$.
Since $T'(Q, n)$ is not homogeneous ideal  
they are also not homogeneous. It is not hard exercise to find the elements
$\tilde{f}$ but it is sufficiently difficult problem to describe homogeneous 
components of their images in $S$. To illustrate this   
take an element $\sigma_r(f)$ from
Theorem 2. For the sake of simplicity we consider only orthogonal invariants. 

Without loss of generality one can assume that $f$ is incident to
the vertex $i_q=1$. If $m$ is a monomial belonging to $f$ then $p_{N, n}(m)=
\tilde{p}_{N, n}(m)$ iff $m$ contains at least one multiplier $Y(a_i), 1\leq
i\leq m$, otherwise $m=(Y(c)Y(b))^l$ or $m=(\overline{Y(c)} \ \overline{Y(b)})
^l$ and $p_{N, n}(m)=\tilde{p}_{N, n}(m) + E(N, n)$. Let $f=f_1+f_2$, where
$f_2$ is a subsum of $f$ which contains all monomials of the second type.
I claim that there are integer coefficients $\alpha_k, 0\leq k\leq r, 
\alpha_0=1,$ such that an element $z=z(f)=\sum_{0\leq k\leq r}\alpha_k
\sigma_{r-k}(f)$ satisfies $p_{N, n}(z)=\sigma_r(\tilde{p}_{N, n}(f))=0$.
Denote by $\lambda$ the sum of all coefficients of the monomials belonging to 
$f_2$. We have 

$$p_{N, n}(t)=\sum_{0\leq j\leq r}\alpha_j\sum_{0\leq k\leq r-j}C^{k}_{N-n}
p_{N, n}(\lambda)^k\sigma_{r-j-k}(\tilde{p}_{N, n}(f))=$$
$$
=\sum_{0\leq t\leq r}\sigma_t(\tilde{p}_{N, n}(f))\sum_{0\leq j, k, j+k=r-t}
C^k_{N-n}p_{N, n}(\lambda)^k\alpha_j .$$

\noindent The required result follows if our coefficients satisfy the 
equations

$$\sum_{0\leq k\leq r-t}
C^k_{N-n}\lambda^k\alpha_{r-t-k}=0, 0 \leq t\leq r-1.$$

\noindent It is clear that these equations has a unique solution 
whenever $\alpha_0$ is fixed, say $\alpha_0=1$. One can prove that
in this case $\alpha_j=(-1)^j C^j_{N-n+j-1}\lambda^j, 0\leq j\leq r$.
For the sake of convenience we denote the $t$-th equation by $Eq_t$.  
\begin{lm}
For any general $N\times N$ matrix $X$ and a variable $y$ we have
$\sigma_k(X+yE(N))=\sum_{0\leq s\leq k}C^s_{N-k+s}y^s\sigma_{k-s}(X)$. 
\end{lm}

Proof. Without loss of generality one can assume that $X$ is a diagonal
matrix with diagonal coefficients $x_1, \ldots , x_N$. Then we have

$$\sigma_k(X+yE(N))=\sum_{1\leq i_1<\ldots < i_k\leq N}(x_{i_1}+y) 
\ldots (x_{i_k}+y)=$$
$$=\sum_{0\leq s\leq k}y^s\sum_{1\leq r_1<\ldots < r_{k-s}\leq
k}x_{i_{r_1}}\ldots x_{i_{r_{k-s}}}.$$

\noindent It is clear that any summand $x_{j_1}\ldots x_{j_{k-s}}, j_1 
<\ldots <  j_{k-s}$, appears as many times as one can choose $s$   
different integers from the set $\{1,\ldots , N\}\setminus\{j_1,\ldots , 
j_{k-s}\}$, that is $C^s_{N-k+s}$. This concludes the proof. 
   
Next, one has to take the matrices $Y(b), Y(c)$ to 
$E(N)$. In particular, $f_2\mapsto\lambda E(N)$. 
We get

$$z'=z(f'_1+\lambda E(N))=\sum_{0\leq k\leq r}\alpha_{r-k}\sum_{0\leq s\leq k}
C^s_{N-k+s}\lambda^s\sigma_{k-s}(f'_1)=$$
$$=\sum_{0\leq t\leq r}\sigma_t(f'_1)
\sum_{0\leq s\leq r-t}C^s_{N-t}\lambda^s\alpha_{r-t-s}.
$$  

\noindent Here $f'_1$ is just the image of $f_1$ under the same specialization 
$Y(b), Y(c)\mapsto E(N)$. 

I claim that all sums $\sum_{0\leq s\leq r-t}C^s_{N-t}\lambda^s\alpha_{r-t-s}$
are equal to zero if $t\leq n$. More generally, 
one can prove that the sums $\sum_{0\leq s
\leq r-t_1}C^s_{N-t_2}\lambda^s\alpha_{r-t_1-s}$ are equal to zero for any pair
$(t_1, t_2)$ such that $0\leq t_1\leq t_2\leq n$. If $t_2=n$ it is just 
the equation $Eq_{t_1}$ (notice that $t_1\leq n\leq r-1$). Let $n > t_2$.
Using the binomial identity $C^j_{n-1}+C^{j-1}_{n-1}=C^j_n$ we have

$$\sum_{0\leq s
\leq r-t_1}C^s_{N-t_2}\lambda^s\alpha_{r-t_1-s}=$$
$$=\sum_{0\leq s
\leq r-t_1}C^s_{N-(t_2+1)}\lambda^s\alpha_{r-t_1-s} +\lambda\sum_{0\leq s
\leq r-(t_1+1)}C^s_{N-(t_2+1)}\lambda^s\alpha_{r-(t_1+1)-s}.$$

\noindent But for pairs $(t_1, t_2+1)$, $(t_1+1, t_2+1)$ 
the induction hypothesis implies that the both last sums are equal to zero.

These computations show that, up to multipliers, all homogeneous components
of $z'$ are $\sigma_j(f'_1), n <j\leq r$.  
It suggets the idea that the defining relations for
the orthogonal or symplectic invariants must be very close to the relations 
from
Theorem 2. But, to realize this idea in a complete form  one has to 
investigate the invariants $\sigma_{r, s}$ more carefully. For example, we 
need some analog of Amitsur's formulae for these invariants. By this reason 
I postpone it for the next article.    
            
\begin{center}
Acknowledgements.
\end{center}

The final part of this research was done during author's vizit to 
Sao Paulo University, Brazil, due to
FAPESP. I am grateful for this support. 
I thank also for RFFI supporting (grant N 01-01-00674). My especially 
thanks to 
Alexandr Grishkov whose efforts played crucial role to make this vizit 
possible.

\end{document}